\documentclass{amsart}

\usepackage[english]{babel}
\usepackage{graphicx}
\usepackage[hidelinks]{hyperref} 
\usepackage{amssymb}
\usepackage{amsmath}
\usepackage[foot]{amsaddr}
\usepackage{amsfonts}
\usepackage[dvipsnames, table]{xcolor}
\usepackage{mathtools}
\usepackage{algpseudocode}
\usepackage{algorithm}
\usepackage{comment}
\usepackage[top=1in, bottom=1.25in, left=1.25in, right=1.25in]{geometry}
\usepackage{float}
\usepackage{tikz}
\newcommand{\norm}[1]{\lVert #1 \rVert}

{}

\newtheorem{remark}{Remark}

\graphicspath{{images/}}

\makeatletter
\algnewcommand{\LineComment}[1]{\Statex \ALG@thistlm {\color{white}\textbf{Input:}} #1}
\makeatother

\makeatletter
\algnewcommand{\LineCommentO}[1]{\Statex \ALG@thistlm {\color{white}\textbf{Output:}} #1}
\makeatother
\algrenewcommand\algorithmicrequire{\textbf{Input:}}
\algrenewcommand\algorithmicensure{\textbf{Output:}}

\begin{document}

\title[FOM and ROM VMS-based EF strategies]{Variational {Multiscale} Evolve and Filter strategies for convection-dominated flows}

\author[M. Strazzullo, F. Ballarin, T. Iliescu, T. Chac\'on Rebollo]{Maria Strazzullo$^{1, a}$, Francesco Ballarin$^{2, a}$, Traian Iliescu$^3$ and  Tom\'as Chac\'on Rebollo$^4$}
\date{\today}
\address{$^1$ Politecnico of Torino, Department of Mathematical Sciences ``Giuseppe Luigi Lagrange'', Corso Duca degli Abruzzi, 24, 10129, Torino, Italy.}
\address{$^2$Department of Mathematics and Physics, Università Cattolica del Sacro Cuore, via Garzetta 48, 25133 Brescia, Italy.}
\address{$^3$Department of Mathematics, Virginia Tech, Blacksburg, VA 24061, USA.}
\address{$^4$ IMUS \& Departamento de Ecuaciones Diferenciales y An\'alisis Num\'erico, Universidad de Sevilla, 41080 Seville, Spain.}
\address{$^{\MakeLowercase{a}}$ INdAM research group GNCS member.}

\begin{abstract}
The evolve-filter (EF) model is a filter-based numerical stabilization for under-resolved
convection-dominated flows. EF is a simple, modular, and effective strategy for both full-order
models (FOMs) and reduced-order models (ROMs). It is well-known, however, that when the
filter radius is too large, EF can be overdiffusive and yield inaccurate results. To alleviate this,
EF is usually supplemented with a relaxation step. The relaxation parameter, however, is very
sensitive with respect to the model parameters. In this paper, we propose a novel strategy to
alleviate the EF overdiffusivity for a large filter radius. Specifically, we leverage the variational
multiscale (VMS) framework to separate the large resolved scales from the small resolved scales in
the evolved velocity, and we use the filtered small scales to correct the large scales. Furthermore, in
the new VMS-EF strategy, we use two different ways to decompose the evolved velocity: the VMS
Evolve-Filter-Filter-Correct (VMS-EFFC) and the VMS Evolve-Postprocess-Filter-Correct (VMS-
EPFC) algorithms. The new VMS-based algorithms yield significantly more accurate results than
the standard EF in both the FOM and the ROM simulations of a flow past a cylinder at Reynolds
number Re = 1000.

\end{abstract}

\maketitle
\section{Introduction}

 The efficient and accurate numerical simulation of the Navier-Stokes equations (NSE) is a central task in many scientific and industrial applications. However, numerical discretizations such as finite element (FE) methods often lead to high computational costs. One of the main challenges is the efficient and accurate simulation of convection-dominated flows, where complex flow structures should be captured down to the Kolmogorov scale. Using such a spatial resolution can yield very large dimensional algebraic systems that can be challenging to solve with the available computational resources. Using coarse meshes in the full-order model (FOM), i.e., working in the \emph{under-resolved} regime, is not a practical option: in this setting, the simulations are usually inaccurate and display spurious numerical oscillations. To address this issue, {numerical stabilization} methods can be employed. For a survey on FOM numerical stabilization approaches, we refer to the monograph \cite{roos2008robust}.

This paper focuses on regularized strategies, which are a particular class of stabilization methods. Their objective is to smooth {out} the numerical oscillations through spatial filtering, while still preserving accuracy in the flow representation. Despite the employment of regularized models on coarse meshes, computational fluid dynamics simulations remain costly and unfeasible in many applied fields. To overcome this issue, \textit{reduced order models} (ROMs) can be used. ROMs are a surrogate version of the FOM and decrease the number of degrees of freedom of the FOM space by orders of magnitude \cite{benner2017model, hesthaven2015certified, prud2002reliable, quarteroni2015reduced, RozzaHuynhPatera2008}. The ROM is a low-dimensional model that exploits the FOM information in the definition of basis functions. Once built, the ROM is used in a Galerkin framework to perform simulations in a {shorter} time than the FOM. The main objective of ROMs is to make an impact in applications, accelerating the simulation process while still being {sufficiently} accurate. In literature it is well established that, in a convection-dominated regime, if stabilization or regularization is needed at the FOM level, it might also be beneficial at the ROM level \cite{girfoglio2021pod, pacciarini2014stabilized, Strazzullo20223148,wang2011two,wang2012proper,zoccolan2, zoccolan1}. Regularized ROMs (Reg-ROMs) are applied as in the FOM case, to be accurate in an under-resolved ROM regime, i.e., when the number of ROM basis functions is not high enough to capture the flow dynamics. 
This is shown, for example, for Leray ROM in \cite{kaneko2020towards,sabetghadam2012alpha,wells2017evolve}, for EFR-ROM in \cite{girfoglio2021pod,Strazzullo20223148,wells2017evolve}, and for the approximate deconvolution Leray ROM in \cite{SanFilippo}.

A comprehensive survey of regularized strategies is the seminal monograph of Layton and Rebholz~\cite{layton2012approximate}. Among the proposed techniques, the evolve and filter (EF) approach is particularly successful thanks to its effectiveness and simplicity, see, e.g., \cite{boyd2001chebyshev, CHQZ88, fischer2001filter, germano1986differential-b,germano1986differential,GirfoglioQuainiRozza2019,girfoglio2021pod,Strazzullo20223148, mullen1999filtering, pasquetti2002comments}. The EF strategy consists of two steps: (i) a standard FOM simulation is performed in an under-resolved regime (evolve), and (ii) a spatial filter is applied to alleviate its spurious oscillations. The action of the filter is related to the filter radius, $\delta$. One of the EF's drawbacks is that the results are very sensitive with respect to $\delta$. For example, when $\delta$ is too large, the EF strategy is overdiffusive and introduces too much numerical diffusion, yielding inaccurate results. To tackle this problem, some works, such as \cite{ bertagna2016deconvolution, boyd1998two, ervin2012numerical, fischer2001filter,M2AN_2014__48_3_765_0,mullen1999filtering}, propose an additional step, relaxation, to alleviate the overdiffusivity and yield accurate results. 

In this paper, we propose a fundamentally different strategy, in the variational multiscale (VMS) framework. The new approach consists of four steps: (i) we evolve the NSE solution in the under-resolved regime, (ii) we exploit the VMS structure and decompose the evolved velocity into resolved large and small scales, (iii) we filter the resolved small scale, and (iv) we correct the resolved small scales with the filtered ones.

We propose two novel strategies based on \emph{how} we decompose the velocity in the second step: the VMS evolve-filter-filter-correct (VMS-EFFC) and 
the evolve-postprocess-filter-correct (VMS-EPFC). In the VMS-EFFC approach, the resolved large scales are captured by the employment of a first differential filter, while in the VMS-EPFC approach, they are captured by a postprocessing 
of the evolved velocity in a smaller space.

These strategies, at the FOM level, 
are compared to standard FE-Galerkin NSE simulations and to standard EF approaches that are overdiffusive (e.g., when the filter radius is too large). Our numerical investigation shows that the new VMS-EFFC and VMS-EPFC strategies yield drastically more accurate results with respect to the standard EF. Moreover, we propose the VMS-EFFC-ROM and the VMS-EPFC-ROM approaches as ROM counterparts of the VMS-filter strategies. In the numerical investigation we performed, the VMS-based filter approaches outperform Galerkin ROMs (G-ROMs) and EF-ROM, i.e., EF applied at the ROM level.\\

{The following} are the main novelties of this paper:
\begin{itemize}
    \item[$\circ$] We propose novel VMS-based filter strategies, increasing the accuracy of FE simulations in under-resolved regimes.
    \item[$\circ$] We extend the VMS-based filter to the ROM level, providing a suitable tool to represent complex flow behaviors in a low-dimensional setting.
\end{itemize}

Both the VMS~\cite{hughes1996space,hughes1998variational,iliescu2013variational,iliescu2014variational,koc2023residual,rebollo2014mathematical,stabile2019reduced} and the EF~\cite{BIL05,GirfoglioQuainiRozza2019,girfoglio2021pod,layton2008numerical, Strazzullo20223148} approaches have been separately studied, both in the FOM and the ROM contexts. We highlight contributions \cite{LAYTON20113183} and \cite{EROGLU2017350}, where a total decoupled VMS-projection stabilization is performed at the FOM and ROM levels, respectively. We emphasize, however, that, in contrast with the previous literature, this is the first time the EF and VMS strategies are combined in \emph{correction} algorithms.

The paper is organized as follows: In Section \ref{sec:EF}, we introduce the incompressible NSE at the continuous and discrete levels and the EF strategy. Section \ref{sec:VMSEF} concerns the definitions of the VMS-EFFC and VMS-EPFC algorithms. In Section \ref{sec:ROM}, we present the ROM approaches that we investigated, i.e., G-ROMs, EF-ROMs, and the VMS-based filters at the ROM level. FOM results are presented in Section \ref{sec:FOM_results}, while Section \ref{sec:ROM_results} focuses on the numerical comparison of the ROM strategies. Both the FOM and the ROM numerical investigations are for a flow past a cylinder at Reynolds number $Re=1000$. Finally, we draw the conclusions in Section \ref{sec:conc}. Furthermore, the interested reader may find some additional experiments in the following appendices:

\begin{itemize}
    \item[$\circ$] Appendix \ref{app:diff}, which investigates the role of the filter radius in the VMS-EFFC;
    \item[$\circ$] Appendix \ref{app:pressure}, which investigates the role of pressure at the FOM level;
    \item[$\circ$] Appendix \ref{app:pressureROM}, which investigates the role of pressure at the ROM level.
\end{itemize}

\section{The FOM and the Evolve-Filter ({EF}) Strategy}
\label{sec:EF}
In this section, we introduce the continuous model of the incompressible Navier-Stokes equations (NSE), together with its discretization. Moreover, we present the EF strategy, which alleviates the spurious numerical oscillations in convection-dominated regimes.
\subsection{The Continuous Model}
This section focuses on the incompressible NSE, i.e., the mathematical formulation that we use in the numerical investigation. Let us consider a spatial domain 
$\Omega \subset \mathbb R^{2}$. The goal is to approximate the velocity $u(x, t) \doteq u\in \mathbb U \doteq L^2((0,T); H^1(\Omega))$ and the pressure $p(x, t) \doteq  p \in \mathbb Q \doteq L^2((0,T); L^2(\Omega))$ fields that solve:
 \begin{equation}
\label{eq:NSE}
\begin{cases}
 \displaystyle u_t  - \nu \Delta u + (u \cdot \nabla) u + \nabla p = 0 & \text{in }  \Omega \times (0, T), \\
\nabla \cdot u = 0  & \text{in }  \Omega \times (0, T), \\
u = u_D & \text{on } \Gamma_D \times (0, T), \\
\displaystyle -p n + \nu \frac{{\partial} u}{\partial n} = 0  & \text{on }  \Gamma_N \times (0, T), \\
u(x, 0) = u_0 & \text{in }  \Omega,
\end{cases}
\end{equation}
for an initial condition $u_0$, and a given boundary condition $u_D$ on $\Gamma_D$, where $\Gamma_D \subset \partial \Omega$ and $\Gamma_N \subset \partial \Omega$ are non-overlapping boundary portions with Dirichlet and ``free flow" conditions such that $\overline{\Gamma_D} \cup \overline{\Gamma_N} = \partial \Omega$. The interested reader may refer to \cite{quarteroni2009numerical} for more details on the NSE formulation.
We denote with $H^1_{u_D}(\Omega)$ the space containing functions in $H^1(\Omega)$ satisfying the boundary condition on $\Gamma_D$.\\ 
The flow regime is determined by the Reynolds number
$Re \doteq {\overline{U}{L}}/{\nu}$, where $\overline{U}$ and ${L}$ are the characteristic velocity and length scales of the system, respectively, while $\nu$ is the kinematic viscosity. In this contribution, we investigate settings where the inertial forces dominate the viscous {forces}, i.e., we deal with convection-dominated regimes. 

\subsection{The Discrete Model}
We define a 
triangulation $\mathcal T$ over $\Omega$. Its simplicial elements are denoted with $K$. We employ the FE discretization. Moreover, we define the space $$\mathbb X^n = \{v \in C^0(\overline{\Omega}) \text{ s.t. } v_{|_K} \in \mathbb P_n(K), \; \forall K \in \mathcal T\},$$ where $\mathbb P_n$ is the space of polynomials of degree at most $n$. The FE spaces for velocity and pressure are $\mathbb U_h \doteq H^1_{u_D}(\Omega) \cap \mathbb X^m$ and $\mathbb Q_h \doteq  L^2(\Omega) \cap \mathbb X^l$, for $m, l$ positive integers. The FE dimensions are $N_h^u$ and $N^p_h$ for velocity and pressure, respectively.
In time, we perform an implicit Euler scheme with timestep $\Delta t$. We consider the time instances $t_n = n\Delta t$ for $n = 0, \dots, N_T$, with final time $T = N_T\Delta t$. The solution at time $t^n$ is denoted with $u_h(t^n) = u^n$ and $p_h(t^n) = p^{n}$, for velocity and pressure, respectively. The final discrete system at time $t^{n+1}$ is: find $(u^{n+1},p^{n+1}) \in \mathbb U_h \times \mathbb Q_h$ such that
\begin{equation*}
\begin{cases}
 \displaystyle \int_{\Omega} \frac{u^{n+1} - u^{n}}{\Delta t}v^{n+1}_h \; dx
	+ a(u^{n+1},  v^{n+1})
	+ \hat{c}(u^{n+1}; u^{n+1}, v^{n+1}) - b(v^{n+1}, p^{n+1})   = 0 & \forall v^{n+1} \in \mathbb U_h,
	 \\[0.5cm]
b( u^{n+1},  q^{n+1} ) = 0 & \forall q^{n+1} \in \mathbb Q_h,
\end{cases}
\end{equation*}
where, for all $u, v \in \mathbb U_h$ and $p \in \mathbb Q_h$, the bilinear forms are defined as
$$ 
    a(u,v) = \int_{\Omega} \nabla u : \nabla v \; dx, \quad b(v,p) = \int_{\Omega} p \nabla \cdot v \; dx,
    $$ 
and for all $u, v, z \in \mathbb U_h$,
$$
\hat{c}(u; v, z)
    = \frac{1}{2} [c(u; v, z) - c(u;z,v)],
$$
is the skew-symmetric approximation of the nonlinear term of the NSE \cite{layton2008introduction,temam2001navier}, defined as
$$
{c}(w;u,v) = \int_{\Omega} (w \cdot \nabla)u \cdot v\; dx.$$

\subsection{{The {EF} Strategy}}
\label{sec:efr-fom}
When the inertial forces dominate the viscous forces, spurious numerical oscillations generally arise in under-resolved or marginally-resolved discretizations. 
The {EF} strategy can be employed to regularize the under-resolved simulations  and alleviate this issue. The {EF} algorithm at $t^{n+1}$ reads \cite{Strazzullo20223148}:
\begin{eqnarray}
         &	\text{(I)}& \text{\emph{ Evolve}:} \quad 
\begin{cases}
        	 \displaystyle \frac{{w}^{n + 1} - u^n}{\Delta t} + (w^{n+1} \cdot \nabla)    w^{n+1} - \nu \Delta w^{n+1} + \nabla p^{n+1} = 0 & \text{in } \Omega , \vspace{1mm}\\
\nabla \cdot w^{n+1} = 0 & \text{in } \Omega , \vspace{1mm}\\
w^{n+1} = u_D^{n+1} & \text{on } \Gamma_D , \vspace{1mm}\\
\displaystyle -p^{n+1} n + \nu \frac{\partial {w^{n+1}}}{\partial n} = 0  & \text{on } \Gamma_N . \\
\end{cases}
            \label{eqn:ef-rom-1}\nonumber \\[0.3cm]
            &	\text{(II)} &\text{\emph{ Filter:}} \quad
\begin{cases} 
        	 -\delta^2 \, \Delta u^{n+1} +  u^{n+1} = w^{n+1}& \text{in } \Omega , \vspace{1mm}\\
u^{n+1} = u^{n+1}_D &\text{on } \Gamma_D , \vspace{1mm}\\
\displaystyle \frac{\partial u^{n+1}}{\partial  n} = 0  & \text{on } \Gamma_N .
\end{cases}\nonumber
\end{eqnarray}

Step (I) computes an intermediate velocity $w^{n+1}$ for the new time instance $t^{n+1}$, evolving the velocity $u^{n}$ at $t^{n}$. Step (II) applies a \emph{differential filter} (DF)  with filtering radius $\delta$: it smooths {out} the oscillations of the intermediate velocity approximation $w^{n+1}$, which generally arise in step (I), eliminating the high frequencies from $w^{n+1}$ by means of an elliptic operator, and obtaining the filtered solution $u^{n+1}$. The filtered solution is considered as the new input to step (I), at the new time step. The EF strategy is summarized in Algorithm \ref{alg:EF}.
\begin{algorithm}[H]
\caption{The {EF} algorithm}
\label{alg:EF}
\begin{algorithmic}[1]
\Require initial condition $u^0$, number of time instances $N_T$, timestep $\Delta t$, filtering radius $\delta$
\Ensure filtered solution $u^n$, for $n \in \{1,\dots, N_T\}$ 
\Statex
\For{$n \in \{0,\dots, N_T-1\}$}
\State{Evolve the NSE ($u^{n} \rightarrow w^{n+1}$)}
\State{Apply a DF of radius $\delta$ to $w^{n+1}$ ($w^{n+1} \rightarrow \overline{w}^{n+1}$) }
\State{$u^{n+1} = \overline{w}^{n+1}$}
\EndFor
\end{algorithmic}
\end{algorithm}

We remark that although {different filters can be used} (see, e.g., \cite{girfoglio2021pod} and the references therein), {here we focus on the simplest and most widespread option.}

\begin{remark}[grad-div stabilization]
\label{remark:graddiv}
At the FOM level, 
the DF does not preserve the incompressibility constraint when dealing with no-slip boundary conditions \cite{ervin2012numerical, layton2008numerical, van2006incompressibility}. To tackle this issue, in \cite{ervin2012numerical, layton2008numerical}, a Stokes DF is used to guarantee the mass conservation of the filtered solution.
A less expensive and easier-to-code solution is the grad-div stabilization. This technique adds the following term to step (II):
\begin{equation}
\gamma_D \nabla(\nabla \cdot \overline{u}^{n+1}),
    \end{equation}
for a penalization parameter $\gamma_D > 0$. We refer the reader to \cite{decaria2018determination,john2017divergence,layton2009accuracy} for an overview of the grad-div stabilization strategy. 
\end{remark}
In the {EF} setting, the choice of the filtering radius is \emph{essential}. However, few guidelines can be found in the literature, and the \emph{optimal filter radius} is still an open problem in many CFD applications. For example, a common choice for $\delta$ is the minimum radius of the mesh, i.e.\ $h_{min}$. Another possibility is represented by Kolmogorov's scale $L \cdot Re^{-\frac{3}{4}}$. We refer the reader to \cite{bertagna2016deconvolution,GirfoglioQuainiRozza2019} for more insights.
We emphasize that, in many practical settings, the action of $\delta$ might result in overdiffusive simulations, see, e.g., \cite{GirfoglioQuainiRozza2019,Strazzullo20223148}. In the next section, we address the following question: Are there {EF}-based strategies that are more reliable and accurate than the classical {EF} for $\delta$ values that result in an overdiffusive behavior?

\section{The FOM Variational Multiscale Evolve-Filter Algorithms}
\label{sec:VMSEF}
This section presents two novel algorithms based on the variational multiscale (VMS) paradigm. 
We assume that the two discrete spaces can be decomposed {into} the following direct sums \cite{rebollo2014mathematical}:
\begin{equation}
    \label{eq:decomposition}
    \mathbb U_h = \overline{\mathbb U}_h \oplus \mathbb U_h' \quad \text{and} \quad   \mathbb Q_h = \overline{\mathbb Q}_h \oplus \mathbb Q_h',
\end{equation}
such that $u_h = \overline{u}_h + {u}_h' \in \mathbb U_h$ and
$p_h = \overline{p}_h + {p}_h' \in \mathbb Q_h$ present a unique decomposition.
In this setting, $\overline{\mathbb U}_h$ and $\overline{\mathbb Q}_h$ represent the spaces of the large scales of the velocity and pressure, while ${\mathbb U}_h'$ and ${\mathbb Q}_h'$ represent the small resolved scales of the velocity and pressure, respectively.
We propose two different strategies:
\begin{itemize}
    \item[$\circ$] the VMS Evolve-Filter-Filter-Correct ({VMS-EFFC}),
    \item[$\circ$] the VMS Evolve-Postprocess-Filter-Correct ({VMS-EPFC}).
\end{itemize}
We postpone the detailed description of the algorithms to the following subsections. However, the main ingredients of both techniques are: (i) evolve the solution, (ii) define the large scales of the evolved velocity, (iii) filter the small scales, and (iv) define a new corrected velocity as the sum of the large scales and the filtered small scales.
In this way, in the case of overdiffusive filtering, we can recover some important features of the NSE simulation and ensure more physical simulations.

\subsection{The {VMS-EFFC} Approach} In this strategy, we propose a double filtering process for the evolved velocity solution, and a final correction step. In particular, we define the large scales of the resolved velocity by applying a DF of radius $\delta_1$. Specifically, we rely on the following steps at time $t^{n+1}$, given an input $u^{n}$:
\begin{eqnarray}
         &	\text{(I)}& \text{\emph{ Evolve}:} \quad 
\begin{cases}
        	 \displaystyle \frac{{w}^{n + 1} - u^n}{\Delta t} + (w^{n+1} \cdot \nabla)    w^{n+1} - \nu \Delta w^{n+1} + \nabla p^{n+1} = 0 & \text{in } \Omega , \vspace{1mm}\\
\nabla \cdot w^{n+1} = 0 & \text{in } \Omega , \vspace{1mm}\\
w^{n+1} = u_D^{n+1} & \text{on } \Gamma_D , \vspace{1mm}\\
\displaystyle -p^{n+1} n + \nu \frac{\partial {w^{n+1}}}{\partial n} = 0  & \text{on } \Gamma_N . \\
\end{cases}
            \label{eqn:ef-rom-1}\nonumber \\[0.3cm]
            &	\text{(II)} &\text{\emph{ Filter (obtain} $\overline{w}^{n+1}$\emph{):}} \quad
\begin{cases} 
        	 -\delta^2_1 \, \Delta \overline{w}^{n+1} +  \overline{w}^{n+1} = w^{n+1}& \text{in } \Omega , \vspace{1mm}\\
\overline{w}^{n+1} = u^{n+1}_D &\text{on } \Gamma_D , \vspace{1mm}\\
\displaystyle \frac{\partial \overline{w}^{n+1}}{\partial  n} = 0  & \text{on } \Gamma_N .
\end{cases}\nonumber
\end{eqnarray}
{After defining} the large scales as the filtered evolved velocity, the small scales are computed as ${w'}^{n+1} = w^{n+1} - \overline{w}^{n+1}.$ At this point, we proceed with the following two steps:
\begin{eqnarray}
            &	\text{(III)} &\text{\emph{ Filter (obtain} $\overline{w'}^{n+1}$\emph{):}} \quad
\begin{cases} 
        	 -\delta^2_2 \, \Delta \overline{w'}^{n+1} +  \overline{w'}^{n+1} = {w'}^{n+1}& \text{in } \Omega , \vspace{1mm}\\
\overline{w'}^{n+1} = 0 &\text{on } \Gamma_D , \vspace{1mm}\\
\displaystyle \frac{\partial \overline{w'}^{n+1}}{\partial  n} = 0  & \text{on } \Gamma_N .
\end{cases}\nonumber
\label{eqn:ef-rom-3}\nonumber \\[0.3cm]
            &	\text{(IV)} &\text{\emph{ Correct:}} \quad
        	u^{n+1} = \overline{w}^{n+1} + \overline{w'}^{n+1}.
\nonumber
\end{eqnarray}

Although steps (I) and (II) are the same as in {EF}, the goal is different. The DF of radius $\delta_1$ is applied to define the large scales $\overline{w}^{n+1}$ of the intermediate velocity $w^{n+1}$. At this point, exploiting the direct sum property \eqref{eq:decomposition}, we define the small scales as ${u'}^{n+1} = w^{n+1} - \overline{w}^{n+1}$. The small scales are usually noisy and contain the high-frequency components of the system. In step (III), we smooth the small scales to make them less noisy by means of another DF of radius $\delta_2$. Finally, in step (IV), the approximation at the new time instance is given by the sum of the large scales and the filtered small scales of the evolved velocity, i.e. $u^{n+1} = \overline{w}^{n+1} + \overline{w'}^{n+1}$. We recall that the divergence of $u^{n+1}$ is relatively small thanks to the grad-div stabilization applied in steps (II) and (III). The VMS-EFFC strategy is summarized in Algorithm \ref{alg:EFFC}.
\begin{remark}[Boundary conditions]
\label{remark:homobc}
    We remark that in step (III) of the {VMS-EFFC} and {VMS-EPFC} approaches, homogeneous boundary conditions are imposed all over $\Gamma_D$. In this way, in the correction step (IV), the new input $u^{n+1}$ satisfies the Dirichlet boundary condition $u_D^{n+1}$ over $\Gamma_D$. 
\end{remark}

\begin{algorithm}[H]
\caption{The {VMS-EFFC} algorithm}
\label{alg:EFFC}
\begin{algorithmic}[1]
\Require initial condition $u^0$, number of time instances $N_T$, timestep $\Delta t$, filter radii $\delta_1$ and $\delta_2$
\Ensure filtered solution $u^n$, for $n \in \{1,\dots, N_T\}$ 
\Statex
\For{$n \in \{0,\dots, N_T-1\}$}
\State{Evolve the NSE ($u^{n} \rightarrow w^{n+1}$)}
\State{Apply a DF of  radius $\delta_1$ to $w^{n+1}$ ($w^{n+1} \rightarrow \overline{w}^{n+1}$) }
\State{Apply a DF of  radius $\delta_2$ to ${w'}^{n+1} = {w}^{n+1} - \overline{w}^{n+1}$ (${w'}^{n+1} \rightarrow \overline{w'}^{n+1}$) }
\State{$u^{n+1} = \overline{w}^{n+1} + \overline{w'}^{n+1}$}
\EndFor
\end{algorithmic}
\end{algorithm}

\begin{remark}[The choice of $\delta_1$ and $\delta_2$]
\label{rem:deltas}
To the best of our knowledge, the {VMS-EFFC} approach is novel. Consequently, there are no clear guidelines on what values should be assigned to the two filter radii $\delta_1$ and $\delta_2$. Our preliminary numerical investigation in Appendix \ref{app:diff} shows that choosing $\delta_1 \neq \delta_2$ does not yield more accurate results than the $\delta_1 = \delta_2$ choice. Thus, in what follows, we use the same filter radius for both steps (II) and (III), i.e., $\delta_1 = \delta_2 = \delta$.
\end{remark}

\begin{remark}[Relation to Approximate Deconvolution]
{Assuming that $\delta_1=\delta_2=\delta$,} \color{black} the output velocity of the {VMS-EFFC} can be interpreted as a result of a van Cittert Approximate Deconvolution (AD) operator of order 1. This AD operator is defined as \cite{layton2012approximate}
 $$
 D_1(u) = 2 \overline{u} - \overline{\overline{u}},
 $$
 {where, $\overline{u}$ is the filtered solution obtained as in step (II), while $\overline{\overline{u}}$ is obtained by applying the differential filter of step (II) with the $\overline{u}$ velocity on the right-hand side.} \color{black}
We note that the VMS-EFFC velocity, {which applies a first differential filter to find $u'$, a second differential filter to find $\overline{u'}$, and then, corrects the large scales, }\color{black} is formally equal to $ D_1(u)$ since 
 \begin{align}
 u = \overline{u} + \overline{u'} =
 \overline{u} + \overline{u - \overline{u}} = 2 \overline{u} - \overline{\overline{u}} = D_1(u).
 \end{align}
We emphasize, however, that this relation holds for $\delta_1=\delta_2$, but generally not for $\delta_1 \neq \delta_2$.
\end{remark}

\subsection{The {VMS-EPFC} Approach} This strategy proposes a different approach to define the small scales of the evolved velocity. In this case, we employ a restriction operator $\Pi_h: \mathbb U_h \rightarrow \overline{\mathbb U}_h$ such that $\Pi_h u_h = \overline{u}_h$, for all $u_h \in \mathbb U_h$. \\
Specifically, we look for $\overline{u}_h \in \overline{\mathbb U}_h$ such that, for all $v \in \overline{\mathbb U}_h$, the following holds: 
\begin{equation}
\label{eq:proj_graddiv}
    (\overline{u}_h, v)_{L^2(\Omega)} + \gamma_P (\nabla \cdot \overline{u}_h, \nabla \cdot v)_{L^2(\Omega)}   = (u_h, v)_{L^2(\Omega)},
\end{equation}
for $\gamma_P>0$. Namely, we are \emph{postprocessing} the evolved velocity by means of an $L^2$-projection enhanced with grad-div stabilization. We recall that the $L^2$-projection is a classical choice in the context of VMS-based strategy related to Smagorinsky-type models, see, e.g., \cite{LAYTON20113183,rebollo2014mathematical} and the references therein. We note, however, that the $L^2$-projection does not satisfy the incompressibility constraint. Thus, to ensure a physically consistent correction step in the VMS-EPFC algorithm, we need the mass conservation to hold for \emph{both} the large scales and the filtered small scales.  For these reasons, we employ grad-div stabilization that, as mentioned in Remark \ref{remark:graddiv}, is a convenient strategy to recover mass conservation properties. The grad-div parameter $\gamma_P$ is chosen by trial and error, since, to our knowledge, we are the first to employ this stabilization in this context. We postpone the discussion on $\gamma_P$ to Section \ref{sec:FOM_results}.

In particular, the VMS-EPFC algorithm at time $t^{n+1}$ reads: 

\begin{eqnarray}
         &	\text{(I)}& \text{\emph{ Evolve}:} \quad 
\begin{cases}
        	 \displaystyle \frac{{w}^{n + 1} - u^n}{\Delta t} + (w^{n+1} \cdot \nabla)    w^{n+1} - \nu \Delta w^{n+1} + \nabla p^{n+1} = 0 & \text{in } \Omega , \vspace{1mm}\\
\nabla \cdot w^{n+1} = 0 & \text{in } \Omega , \vspace{1mm}\\
w^{n+1} = u_D^{n+1} & \text{on } \Gamma_D , \vspace{1mm}\\
\displaystyle -p^{n+1} n + \nu \frac{\partial {w^{n+1}}}{\partial n} = 0  & \text{on } \Gamma_N . \\
\end{cases}
            \label{eqn:ef-rom-1}\nonumber \\[0.3cm]
            &	\text{(II)} &\text{\emph{ Postprocess (obtain} $\overline{w}^{n+1}$\emph{):}} \quad \Pi_hw^{n+1} = \overline{w}^{n+1}.
\nonumber
\end{eqnarray}
Once the large scales of the evolved velocity are defined, we can exploit the direct sum \eqref{eq:decomposition} to obtain ${w'}^{n+1} = w^{n+1} - \overline{w}^{n+1}.$ We can now apply the same steps as in the VMS-EFFC algorithm. We repeat them for the sake of clarity:
\begin{eqnarray}
            &	\text{(III)} &\text{\emph{ Filter (obtain} $\overline{w'}^{n+1}$\emph{):}} \quad
\begin{cases} 
        	 -\delta^2 \, \Delta \overline{w'}^{n+1} +  \overline{w'}^{n+1} = {w'}^{n+1}& \text{in } \Omega , \vspace{1mm}\\
\overline{w'}^{n+1} = 0 &\text{on } \Gamma_D , \vspace{1mm}\\
\displaystyle \frac{\partial \overline{w'}^{n+1}}{\partial  n} = 0  & \text{on } \Gamma_N .
\end{cases}\nonumber
\label{eqn:ef-rom-3}\nonumber \\[0.3cm]
            &	\text{(IV)} &\text{\emph{ Correct:}} \quad
        	u^{n+1} = \overline{w}^{n+1} + \overline{w'}^{n+1}.
\nonumber
\end{eqnarray}

Step (I) evolves the velocity, obtaining $w^{n+1}$ from an input $u^{n}$. In step (II), we obtain the large scales by solving \eqref{eq:proj_graddiv}.  As a consequence, we can define the small scales $w'^{n+1}$. Step (III) applies a DF of radius $\delta$ to the small scales of the evolved velocity. In this way, the noisy behavior of the small scales is alleviated, and important features of the flow are recovered. Finally, in step (IV) we correct the large scales with the filtered small scales to obtain the new input $u^{n+1}$. We recall that \emph{both} $\overline w^{n+1}$ and $\overline {w'}^{n+1}$ feature small divergence values due to the grad-div stabilization applied in steps (II) and (III). The VMS-EPFC strategy is described in Algorithm \ref{alg:EPFC}.
\begin{algorithm}[H]
\caption{The {VMS-EPFC} algorithm}
\label{alg:EPFC}
\begin{algorithmic}[1]
\Require initial condition $u^0$, number of time instances $N_T$, timestep $\Delta t$, filter radius $\delta$
\Ensure filtered solution $u^n$, for $n \in \{1,\dots, N_T\}$ 

\Statex
\For{$n \in \{0,\dots, N_T-1\}$}
\State{Evolve the NSE ($u^{n} \rightarrow w^{n+1}$)}
\State{Apply postprocessing \eqref{eq:proj_graddiv} to $w^{n+1}$ ($w^{n+1} \rightarrow \overline{w}^{n+1}$) }
\State{Apply a DF of  radius $\delta$ to ${w'}^{n+1} = {w}^{n+1} - \overline{w}^{n+1}$ (${w'}^{n+1} \rightarrow \overline{w'}^{n+1}$) }
\State{$u^{n+1} = \overline{w}^{n+1} + \overline{w'}^{n+1}$}
\EndFor
\end{algorithmic}
\end{algorithm}
\section{The ROM Evolve-Filter and Variational Multiscale Evolve-Filter Algorithms}
\label{sec:ROM}
In this section, we extend the VMS strategies that we proposed in Section \ref{sec:VMSEF} to reduced order modeling. We first present the techniques commonly used in the literature, i.e., Galerkin-ROM (G-ROM) and EF-ROM. Then, we propose the ROM versions of the variational multiscale strategies of Section \ref{sec:VMSEF}: VMS-EFFC-ROM and VMS-EPFC-ROM.

\subsection{The G-ROM Approach}
\label{sec:g-rom}
This section presents the G-ROM strategy. At the FOM level, a regularized formulation is used. Specifically, we obtain the FOM solutions applying the VMS-EFFC regularization. At the reduced level, a standard Galerkin projection over a ROM space is performed in the online phase. We exploit the proper orthogonal decomposition (POD) algorithm, which builds the reduced spaces by means of compression of FOM data for several time instances, i.e., the \emph{snapshots}. 
The POD yields two sets of basis functions $\{\phi_j\}_{j=1}^{r_u}$ and $\{\psi_j\}_{j=1}^{r{_p}}$, spanning the reduced spaces $\mathbb U^{r_u}$ and $\mathbb Q^{r_p}$, respectively. 
The reduced velocity $u_{r_u} \in \mathbb U^{r_u}$ and the reduced pressure $p_{r_p} \in \mathbb Q^{r_p}$ can be expanded as
$ 
	u_{r_u} \doteq {u}^{r_u}(x,t)
	= \sum_{j=1}^r a_j^{u}(t) \phi_j(x)$ {and} $
  p_{r_p} \doteq {p}^{r_p}(x,t)
	= \sum_{j=1}^r a_j^{p}(t) \psi_j(x),
	\label{eqn:g-rom-1}
$ 
for some reduced real coefficients $\{a_{j}^{u}(t)\}_{j=1}^{r_u}$ and $\{a_{j}^{p}(t)\}_{j=1}^{r_p}$, respectively \cite{noack2011reduced}. We build the spaces considering $N_{u}$ and $N_p$ snapshots for velocity and pressure, respectively, i.e.,
$\{u_i\}_{i=1}^{N_{u}} \subseteq \{u^k\}_{k=1}^{N_T}$ and $\{p_i\}_{i=1}^{N_{p}} \subseteq \{p^k\}_{k=1}^{N_T}$. {In our numerical tests, we consider a subset made of equispaced snapshots in time}.
The POD algorithm \cite{burkardt2006pod, hesthaven2015certified}  is enhanced by \emph{supremizer stabilization}, which enriches the reduced velocity space $\mathbb U^r$ to guarantee the well-posedness of the projected system \cite{ballarin2015supremizer,rozza2007stability}, avoiding spurious reduced pressure modes.
However, in a convection-dominated setting, other stabilization approaches are generally needed at the reduced level. This will be the topic of Sections \ref{sec:ef}, \ref{sec:VMS-EFFC-ROM}, and \ref{sec:VMS-EPFC-ROM}.
The supremizer operator is defined as
$S: \mathbb Q_h \rightarrow{{\mathbb U}}_h$ with
$(S(p), \tau)_{\mathbb {U}} = \left ( p, \nabla \cdot \tau \right )$ for all $\tau \in \mathbb{U}_h$. 
The final reduced velocity space is
$
{{\mathbb U}}^{r_{us}} \doteq \text{POD}(\{u_i\}_{i=1}^{N_{u}}, r_u)\oplus \text{POD}(\{S(p_i)\}_{i=1}^{N_{u}}, r_s).
$
The supremizer operator is computed for each pressure related to the velocity snapshots.
The reduced pressure is obtained by the standard POD algorithm
$
\mathbb Q^{r_p} \doteq \text{POD}(\{p_i\}_{i=1}^{N_{p}}, r_p).
$
In the compression process, we retain $r_u$, $r_s$, and $r_p$ basis functions for the velocity, supremizer, and pressure, respectively. The final reduced velocity space is spanned by $r_{us} = r_u + r_s$ modes.
Thus, after performing the POD and the supremizer enrichment, we employ the bases $\{\phi_j\}_{j=1}^{r_{us}}$ 
and $\{\psi_j\}_{j=1}^{r_p}$ for the velocity and pressure, respectively.
After building the velocity and pressure bases, at each time instance $t^{n+1}$, we solve the evolution of the system projecting it on the reduced spaces: 
\begin{equation}
\begin{cases}
	\displaystyle \int_{\Omega}\frac{u_{r_{us}}^{n+1} - u_{r_{us}}^{n}} {\Delta t} \phi_i \; dx
	+ \nu a(u_{r_{us}}^{n+1},  \phi_i)
	+ \hat{c}(u_{r_{us}}^{n+1}; u_{r_{us}}^{n+1} , \phi_i)  - b(\phi_i, p_{r_p}^{n+1})
	= 0, \\[0.5cm]
b( u_{r_{us}}^{n+1},  \psi_j ) = 0,
\end{cases}
\label{eqn:NSE_r}
\end{equation}
for all {$\phi_i$ with $ i = 1, \ldots, r_{us},$ and $\psi_j$ with $j = 1, \ldots, r_p$}, where the G-ROM solution is expanded as
\begin{equation}
	{u}_{r_{us}} = \sum_{i=1}^{r_{us}} a_i^{{u}}(t) \phi_i(x) \quad \text{and} 
 \quad
 {p}_{r_{p}} = \sum_{j=1}^{r_{p}} a_i^{{p}}(t) \psi_i(x).
	\label{eqn:us}
\end{equation}
The G-ROM strategy is summarized in Algorithm \ref{alg:GROM}.
\begin{algorithm}[H]
\caption{The {G-ROM} algorithm}
\label{alg:GROM}
\begin{algorithmic}[1]
\Require initial condition $u^0$, number of time instances $N_T$, timestep $\Delta t$, filter radius $\delta$, number of snapshots
\Ensure reduced solution $u^n_{r_{us}}$, for $n \in \{1,\dots, N_T\}$ 
\Statex 
\State apply Algorithm \ref{alg:EFFC}
\State{$\{u_i\}_{i=1}^{N_{u}} \subseteq \{u^k\}_{k=1}^{N_{T}}$ 
\quad $\{p_i\}_{i=1}^{N_{p}} \subseteq \{p^k\}_{k=1}^{N_{T}}$} \State{$\mathbb U^{r_{us}} \doteq \text{POD}(\{u_i\}_{i=1}^{N_{u}})\oplus \text{POD}(\{S(p_i)_{i=1}^{N_{u}})\}$ } \State{$\mathbb Q^{r_p} \doteq \text{POD}(\{p_i\}_{i=1}^{N_p})$ } \For{$n \in \{0,\dots, N_T-1\}$}
\State{Solve \eqref{eqn:NSE_r} }
\EndFor
\end{algorithmic}
\end{algorithm}

\subsection{The EF-ROM approach}
	\label{sec:ef}
The EF-ROM strategy employs the EF algorithm at the ROM level. The reduced space is built as shown in Section \ref{sec:g-rom}, applying the POD and the supremizer enrichment on the VMS-EFFC solutions. At the online stage, we first evolve in the reduced {space} and, then, filter the solution. Namely, we solve:

\begin{eqnarray*}
         &	\text{(I)}_r& \text{\emph{}} \quad 
\begin{cases}
	\displaystyle \int_{\Omega}\frac{{w}_{r_{us}}^{n+1} - u_{r_{us}}^{n}} {\Delta t} \phi_i \; dx
	+ \nu a({{w}}_{r_{us}}^{n+1},  \phi_i)
	+ \hat{c}({{w}}_{r_{us}}^{n+1}; {{w}}_{r_{us}}^{n+1} , \phi_i)  - b(\phi_i, p_{r_p}^{n+1})
	= 0, \\[0.5cm]
b( {{w}}_{r_{us}}^{n+1},  \psi_j ) = 0,
\end{cases}
            \label{eqn:ef-rom-1-r}\nonumber \\[0.3cm]
            &	\text{(II)}_r &\text{\emph{}} \quad
        	 \delta^2 \, a(\overline{w}^{n+1}_{r_{us}}, \phi_i)  +  \int_{\Omega} \overline{w}^{n+1}_{r_{us}} \phi_i \; dx = \int_{\Omega} {w}^{n+1}_{r_{us}} \phi_i \; dx, 
        	\\
	\label{eqn:ef-rom-2-r} \nonumber 
\end{eqnarray*}
 for all {$\phi_i$ with $ i = 1, \ldots, r_{us},$ and $\psi_j$ with $j = 1, \ldots, r_p$}. At each time instance, we expand the reduced variables $\overline{w}_{r_{us}}^n$ and ${w}_{r_{us}}^n$ in $\mathbb U^{r_{us}}$ as \begin{equation}
	\overline{w}_{r_{us}}^n
	= \sum_{j=1}^{r_{us}} a_j^{\overline{w}}(t) \phi_j(x) \quad \text{and} \quad {w}^n_{r_{us}} 
	= \sum_{j=1}^{r_{us}} a_j^{{w}}(t) \phi_j(x),
	\label{eqn:ws}
\end{equation}
and we set the next reduced input as $u_{r_{us}}^{n+1} = \overline{w}_{r_{us}}^{n+1}$. We summarize the EF-ROM strategy in Algorithm \ref{alg:EF-ROM}. We remark that in this contribution we used the same $\delta$ for FOM and ROM simulations, see \cite{Strazzullo20223148}. {Yet another choice is based on energy balancing arguments~\cite{mou2023energy}.}
\begin{algorithm}[H]
\caption{The {EF-ROM} algorithm}
\label{alg:EF-ROM}
\begin{algorithmic}[1]
\Require initial condition $u^0$, number of time instances $N_T$, timestep $\Delta t$, filter radius $\delta$, number of snapshots
\Ensure reduced solution $u^n_{r_{us}}$, for $n \in \{1,\dots, N_T\}$ 
\Statex 
\State apply Algorithm \ref{alg:EFFC}
\State{$\{u_i\}_{i=1}^{N_{u}} \subseteq \{u^k\}_{k=1}^{N_{T}}$ 
\quad $\{p_i\}_{i=1}^{N_{p}} \subseteq \{p^k\}_{k=1}^{N_{T}}$} \State{$\mathbb U^{r_{us}} \doteq \text{POD}(\{u_i\}_{i=1}^{N_{u}})\oplus \text{POD}(\{S(p_i)_{i=1}^{N_{u}})\}$ } \State{$\mathbb Q^{r_p} \doteq \text{POD}(\{p_i\}_{i=1}^{N_p})$ } \For{$n \in \{0,\dots, N_T-1\}$}
\State{Evolve the reduced NSE ($u^{n}_{r_{us}} \rightarrow w^{n+1}_{r_{us}}$)}
\State{Apply a DF of radius $\delta$ to $w^{n+1}_{r_{us}}$ ($w^{n+1}_{r_{us}} \rightarrow \overline{w}^{n+1}_{r_{us}}$)}
\State{$u^{n+1}_{r_{us}} = \overline{w}^{n+1}_{r_{us}}$}
\EndFor
\end{algorithmic}
\end{algorithm}

\subsection{The VMS-EFFC-ROM Approach}
\label{sec:VMS-EFFC-ROM}
In this section, we extend the VMS-EFFC strategy to the reduced level (VMS-EFFC-ROM). To do so, we decompose the reduced velocity space in the direct sum 
\begin{equation}
    \label{eq:decomposition_rom}
    \mathbb U^{r_{us}} = \overline{\mathbb{ U}^{r_{us}}} \oplus \mathbb U^{{r_{us}'}},
\end{equation}
such that the reduced velocity field can be decomposed as $u_{r_{us}} = \overline{u}_{r_{us}} + {u'}_{r_{us}} \in \mathbb U^{r_{us}}$. In this setting, $\overline{\mathbb{ U}^{r_{us}}}$ is the reduced space modeling the large scales of the flow, while $\mathbb U^{{r_{us}'}}$ models the small resolved scales for the reduced velocity.
In this approach, the large scales are derived by applying a DF of radius $\delta_1$. After performing the POD on a FOM discretization, we apply the following steps at time $t^{n+1}$, given a reduced input $u_{r_{us}}^{n}$:

\begin{eqnarray*}
         &	\text{(I)}_r& \text{\emph{}} \quad 
\begin{cases}
	\displaystyle \int_{\Omega}\frac{{w}_{r_{us}}^{n+1} - u_{r_{us}}^{n}} {\Delta t} \cdot \phi_i \; dx
	+ \nu a({{w}}_{r_{us}}^{n+1},  \phi_i)
	+ \hat{c}({{w}}_{r_{us}}^{n+1}; {{w}}_{r_{us}}^{n+1} , \phi_i)  - b(\phi_i, p_{r_p}^{n+1})
	= 0, \\[0.5cm]
b( {{w}}_{r_{us}}^{n+1},  \psi_j ) = 0,
\end{cases}
            \label{eqn:ef-rom-1-r}\nonumber \\[0.3cm]
            &	\text{(II)}_r &\text{\emph{}} \quad
        	 \delta^2_1 \, a(\overline{w}^{n+1}_{r_{us}}, \phi_i)  +  \int_{\Omega} \overline{w}^{n+1}_{r_{us}} \cdot \phi_i \; dx = \int_{\Omega} {w}^{n+1}_{r_{us}} \cdot \phi_i \; dx,
        	\\
	\label{eqn:ef-rom-2-r} \nonumber 
\end{eqnarray*}
for all {$\phi_i$ with $ i = 1, \ldots, r_{us},$ and $\psi_j$ with $j = 1, \ldots, r_p$}. 
In other words, we evolve the reduced velocity and we filter it by means of a DF of radius $\delta_1$. Once the large scales are obtained, we define ${w'}^{n+1}_{r_{us}} = {w}^{n+1}_{r_{us}} - \overline{w}^{n+1}_{r_{us}}$ exploiting decomposition \eqref{eq:decomposition_rom}. 
The VMS-EFFC-ROM strategy continues with the following two steps:
\begin{eqnarray} \nonumber
            &	\text{(III)}_r &
        	 \delta^2_2 \, a(\overline{w'}^{n+1}_{r_{us}}, \phi_i)  +  \int_{\Omega} \overline{w'}^{n+1}_{r_{us}} \cdot \phi_i \; dx = \int_{\Omega} {w'}^{n+1}_{r_{us}} \cdot \phi_i \; dx,
\end{eqnarray}
for all {$\phi_i$ with $ i = 1, \ldots, r_{us},$ and $\psi_j$ with $j = 1, \ldots, r_p$}, and
\begin{eqnarray} \nonumber
            &	\text{(IV)}_r & \quad
        	u^{n+1}_{r_{us}} = \overline{w}^{n+1}_{r_{us}} + \overline{w'}^{n+1}_{r_{us}}.
\end{eqnarray}

The VMS-EFFC-ROM strategy is summarized in Algorithm \ref{alg:EFFC-ROM}.
\begin{algorithm}[H]
\caption{The {VMS-EFFC-ROM} algorithm}
\label{alg:EFFC-ROM}
\begin{algorithmic}[1]
\Require initial condition $u^0$, number of time instances $N_T$, timestep $\Delta t$, filter radii $\delta_1$ and $\delta_2$, number of snapshots
\Ensure reduced solution $u^n_{r_{us}}$, for $n \in \{1,\dots, N_T\}$ 
\Statex 
\State Apply Algorithm \ref{alg:EFFC}
\State{$\{u_i\}_{i=1}^{N_{u}} \subseteq \{u^k\}_{k=1}^{N_{T}}$ 
\quad $\{p_i\}_{i=1}^{N_{p}} \subseteq \{p^k\}_{k=1}^{N_{T}}$} 
\State{$\mathbb U^{r_{us}} \doteq \text{POD}(\{u_i\}_{i=1}^{N_{u}})\oplus \text{POD}(\{S(p_i)_{i=1}^{N_{u}})\}$ } \State{$\mathbb Q^{r_p} \doteq \text{POD}(\{p_i\}_{i=1}^{N_p})$ }
\For{$n \in \{0,\dots, N_T-1\}$}
\State{Evolve the reduced velocity ($u^{n}_{r_{us}} \rightarrow w^{n+1}_{r_{us}}$)}
\State{Apply a DF of radius $\delta_1$ to $w^{n+1}_{r_{us}}$ ($w^{n+1}_{r_{us}} \rightarrow \overline{w}^{n+1}_{r_{us}}$) }
\State{Apply a DF of radius $\delta_2$ to ${w'}^{n+1}_{r_{us}} = {w}^{n+1}_{r_{us}} - \overline{w}^{n+1}_{r_{us}}$ (${w'}^{n+1}_{r_{us}} \rightarrow \overline{w'}^{n+1}_{r_{us}}$) }
\State{$u^{n+1}_{r_{us}} = \overline{w}^{n+1}_{r_{us}} + \overline{w'}^{n+1}_{r_{us}}$}
\EndFor
\end{algorithmic}
\end{algorithm}
\begin{remark}[Choosing the reduced filter radii]
As mentioned in Remark \ref{rem:deltas}, there are no clear guidelines for choosing the two filter radii $\delta_1$ and $\delta_2$ at the FOM level. At the ROM level, our preliminary numerical investigation shows that choosing $\delta_1 \neq \delta_2$ does not yield accurate results (just as the FOM setting). Thus, for the sake of brevity, in Section \ref{sec:ROM_results}, we do not show the results. This inaccurate behavior for different filter radii is not unexpected, see, e.g., \cite{Strazzullo20223148,zoccolan2,zoccolan1}. Guided by these works, we used a \emph{consistent VMS-EFFC-ROM}, choosing the same FOM filter radii, i.e., $\delta_1 = \delta_2 = \delta$. 
\end{remark}
\subsection{The VMS-EPFC-ROM Approach}
\label{sec:VMS-EPFC-ROM}
This section extends the VMS-EPFC strategy to the ROM setting, i.e., it develops the VMS-EPFC-ROM algorithm. In this case, we use a ROM projection to define the large scale of the evolved velocity. We note that this projection preserves the weakly divergence-free property of the underlying snapshots. Thus, we do not need to use the grad-div stabilization, as we did in the VMS-EPFC algorithm at the FOM level. We proceed by choosing $\overline r_{u} \leq r_{u}$ and defining $r'_{u} = r_{u} - \overline r_{u}$. Recall that, due to the supremizer enrichment, expansion \eqref{eqn:ws} becomes
\begin{equation}
\begin{split}
    	\overline{w}_{r_{us}}^n
	= \sum_{j=1}^{r_{u}} a_j^{\overline{w}_u}(t) \phi_j(x) + \sum_{j=r_u + 1}^{r_{us}} a_k^{\overline{w}_s}(t) \phi_k(x), \\ 
    {w}^n_{r_{us}} 
	= \sum_{j=1}^{r_{u}} a_j^{{w}_u}(t) \phi_j(x) + \sum_{j=r_u + 1}^{r_{us}} a_k^{{w}_s}(t) \phi_k(x),
	\label{eqn:ws_sup}
    \end{split}
\end{equation}
where the last modes $\{\phi_k\}_{k=r_u + 1}^{r_{us}}$ are the supremizer modes. We note that the supremizers include the small scales of the velocity. We refer the reader to Remark \ref{rem:role_proj} for more details about this choice in our numerical investigation. 
The idea of the postprocess step in the VMS-EPFC-ROM is to discard the last $r'_u$ modes of the evolved velocity to define the large scales.
Namely, expanding the evolved velocity as in \eqref{eqn:ws_sup}, we take the following two steps, for all {$\phi_i$ with $ i = 1, \ldots, r_{us},$ and $\psi_j$ with $j = 1, \ldots, r_p$}:

\begin{eqnarray*}
         &	\text{(I)}_r& \text{\emph{}} \quad 
\begin{cases}
	\displaystyle \int_{\Omega}\frac{{w}_{r_{us}}^{n+1} - u_{r_{us}}^{n}} {\Delta t} \phi_i \; dx
	+ \nu a({{w}}_{r_{us}}^{n+1},  \phi_i)
	+ \hat{c}({{w}}_{r_{us}}^{n+1}; {{w}}_{r_{us}}^{n+1} , \phi_i)  - b(\phi_i, p_{r_p}^{n+1})
	= 0, \\[0.5cm]
b( {{w}}_{r_{us}}^{n+1},  \psi_j ) = 0,
\end{cases}
            \label{eqn:ef-rom-1-r}\nonumber \\[0.3cm]
            &	\text{(II)}_r &  
            a_j^{\overline{w}_u}(t) = a_j^{w}(t) \quad \text{for} \quad j= 1, \dots, \overline{r}_{u}, \quad \text{and} \quad 
            a_j^{\overline{w}_u}(t) = 0 \quad \text{for} \quad j= \overline{r}_{u} + 1, \dots, r_{u}.
        	\\
	\label{eqn:ef-rom-2-r} \nonumber 
\end{eqnarray*}
After step (II)$_r$, we define the large scales and the small scales, respectively, as 

\begin{equation}
\begin{split}
    	\overline{w}^{n+1}_{r_{us}} & = \sum_{j=1}^{\overline{r}_{u}} a_j^{{w}_u}(t) \phi_j(x)  + \sum_{j=r_u + 1}^{r_{us}} a_k^{\overline{w}_s}(t) \phi_k(x) \\
 {w'}^{n+1}_{r_{us}} & = \sum_{j=\overline{r}_{u} + 1}^{{r}_{u}} a_j^{{w}_u}(t) \phi_j(x).
	\label{eqn:us}
    \end{split}
\end{equation}
Next, we perform steps (III)$_r$ with a filter radius $\delta$ and (IV)$_r$, as already described in Section \ref{sec:VMS-EFFC-ROM} for VMS-EFFC-ROM, to obtain a new input for the evolve step. The VMS-EPFC-ROM approach is summarized in Algorithm \ref{alg:EPFC-ROM}.
\begin{table}[H]
\caption{Acronyms of the FOM and ROM strategies. Gray cells 
in the ROM-column denote that no ROM is performed.}
\vspace{2mm}
\label{tab:Acro}
\resizebox{\textwidth}{!}{
\begin{tabular}{|c|c|c|c|c|}
\hline
Acronym & FOM regularization&  ROM regularization & VMS regularization and correction & Algorithm\\ \hline
EF &  $\checkmark$ &\cellcolor{gray!50} & &\ref{alg:EF}\\ \hline
VMS-EFFC &  $\checkmark$ & \cellcolor{gray!50} &$\checkmark$  &  \ref{alg:EFFC}\\ \hline
VMS-EPFC  &  $\checkmark$ & \cellcolor{gray!50} & $\checkmark$ &  \ref{alg:EPFC} \\ \hline
G-ROM  &  $\checkmark$ & $\checkmark$ & &  \ref{alg:GROM} \\ \hline
EF-ROM  &  $\checkmark$ & $\checkmark$ &  &  \ref{alg:EF-ROM} \\ \hline
VMS-EFFC-ROM  &  $\checkmark$ & $\checkmark$ & $\checkmark$ & \ref{alg:EFFC-ROM} \\ \hline
VMS-EPFC-ROM  &  $\checkmark$ & $\checkmark$ & $\checkmark$ &  \ref{alg:EPFC-ROM} \\ \hline

\end{tabular}}
\end{table}

\begin{algorithm}[H]
\caption{The {VMS-EPFC-ROM} algorithm}
\label{alg:EPFC-ROM}
\begin{algorithmic}[1]
\Require initial condition $u^0$, number of time instances $N_T$, timestep $\Delta t$, filter radius $\delta$, number of snapshots
\Ensure reduced solution $u^n_{r_{us}}$, for $n \in \{1,\dots, N_T\}$ 
\Statex 
\State apply Algorithm \ref{alg:EFFC}
\State{$\{u_i\}_{i=1}^{N_{u}} \subseteq \{u^k\}_{k=1}^{N_{T}}$ 
\quad $\{p_i\}_{i=1}^{N_{p}} \subseteq \{p^k\}_{k=1}^{N_{T}}$}
\State{$\mathbb U^{r_{us}} \doteq \text{POD}(\{u_i\}_{i=1}^{N_{u}})\oplus \text{POD}(\{S(p_i)_{i=1}^{N_{u}})\}$ } \State{$\mathbb Q^{r_p} \doteq \text{POD}(\{p_i\}_{i=1}^{N_p})$ }
\For{$n \in \{0,\dots, N_T-1\}$}
\State{Evolve the reduced velocity ($u^{n}_{r_{us}} \rightarrow w^{n+1}_{r_{us}}$)}
\State{Set $a_j^{\overline{w}}(t) = a_j^{w}(t) \text{ for } j= 1, \dots, \overline{r}_{us}\text{ and }
            a_j^{\overline{w}}(t) = 0 \text{ for } j= \overline{r}_{us} + 1, \dots, r_{us},$ }
\State{Apply a DF of radius $\delta$ to ${w'}^{n+1}_{r_{us}} = {w}^{n+1}_{r_{us}} - \overline{w}^{n+1}_{r_{us}}$ (${w'}^{n+1}_{r_{us}} \rightarrow \overline{w'}^{n+1}_{r_{us}}$) }
\State{$u^{n+1}_{r_{us}} = \overline{w}^{n+1}_{r_{us}} + \overline{w'}^{n+1}_{r_{us}}$}
\EndFor
\end{algorithmic}
\end{algorithm}

The acronyms and features of the FOM and ROM algorithms are summarized in Table \ref{tab:Acro}.

\section{FOM Results}
\label{sec:FOM_results}
In this section, we perform a numerical investigation to understand which strategy performs better among EF, VMS-EFFC, and VMS-EPFC in the context of convection-dominated NSE.\\
We first describe the general setting and, then, compare the three approaches in terms of accuracy of the velocity fields with respect to a direct numerical simulation (DNS) of the NSE velocity. \\
Additional results are presented in Appendices \ref{app:diff} and \ref{app:pressure}, where we test the VMS-EFFC approach for $\delta_1 \neq \delta_2$, and present results for the pressure fields, respectively.

\subsection{General Numerical Setting and Motivation For Filter Use}
\label{sec:gen_setting}
We introduce the computational setting we investigate in the numerical tests.
We solve NSE \eqref{eq:NSE} on the following spatial domain: $\Omega \doteq \{(0, 2.2) \times (0, 0.41)\} \setminus \{(x, y)\in \mathbb R^2 \text{ such that } (x - 0.2)^2 + (y - 0.2)^2 - 0.05^2 \leq 0 \}$, which is depicted in Figure \ref{fig:domain}. 
The problem features the following Dirichlet boundary condition: 
\begin{equation}
\label{eq:inlet}
u_D = 
\begin{cases}
0 & \text{ on } \Gamma_{\text{w}}, \\
u_{\text{in}} = \displaystyle \left  ( \frac{6}{0.41^2}y(0.41 - y), 0 \right) & \text{ on } \Gamma_{\text{in}}, 
\end{cases}
\end{equation}
where $\Gamma_{\text{w}}$ (solid cyan line in Figure \ref{fig:domain}) is the union of the bottom, 
top, and cylinder walls. The non-homogeneous inlet  condition is applied to $\Gamma_{\text{in}} = \{0\} \times [0, 0.41]$, (dotted magenta line in Figure \ref{fig:domain}), as prescribed in equation \eqref{eq:inlet}.
Moreover, $\Gamma_N$ (dashed black line in Figure \ref{fig:domain}) features
{``free flow"} boundary conditions. 
We use the initial condition  $u_0 = (0,0)$.

\begin{figure}[H]
\begin{center}
\begin{tikzpicture}[scale=5.0]

\filldraw[color=cyan!80, fill=gray!10, very thick](0.2,0.2) circle (0.05);
\filldraw[color=magenta!90, very thick, dotted](0,0) -- (0.,0.41);
\filldraw[color=cyan!80, fill=gray!10, very thick](0,0.41) -- (2.2,0.41);
\filldraw[color=cyan!80, fill=gray!10, very thick](0,0.) -- (2.2,0.);
\filldraw[color=black!80, fill=gray!10, very thick, dashed](2.2,0.) -- (2.2,0.41);

\node at (-.1,0.205){\color{black}{$\Gamma_{\text{in}}$}};
\node at (2.3,0.205){\color{black}{$\Gamma_{N}$}};
\node at (0,-.05){\color{black}{$(0,0)$}};

\node at (0,0.45){\color{black}{$(0,0.41)$}};
\node at (2.35,0.45){\color{black}{$(2.2,0.41)$}};
\node at (2.3,-0.05){\color{black}{$(2.2,0)$}};
\node at (1.1,-.05){\color{black}{$ \color{cyan}{\Gamma_\text{w}}$}};
\end{tikzpicture}
\end{center}
\caption{Spatial domain $\Omega$: schematic representation. $\Gamma_D = \Gamma_{\text{in}} \cup \Gamma_\text{w}$. Homogeneous Dirichlet conditions are applied on the solid cyan boundary. The inlet boundary $\Gamma_{\text{in}}$ is represented by a dotted magenta line. The 
{``free flow"} boundary $\Gamma_N$ is depicted by a dashed black line.}
\label{fig:domain}
\end{figure}
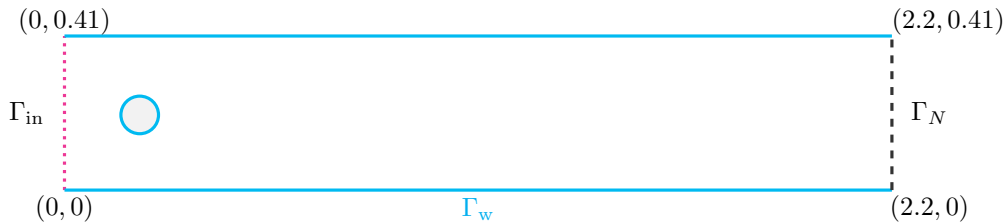
 \begin{table}[H]
\caption{Ganeral setting. Comparison of the coarse mesh (M1), and the fine mesh (M2) with respect to $h_{min}$, $h_{max}$, and number of 
{degrees of freedom}.}
\vspace{3mm}
\label{tab:meshval} \centering
\begin{tabular}{|c|c|c|c|}
\hline
mesh & $h_{min}$            & $h_{max}$            &  $N_h$     \\ \hline
M1   & $4.46\cdot 10^{-3}$  & $4.02 \cdot 10^{-2}$ & $14053$  \\ \hline
M2   & $4.44 \cdot 10^{-3}$ & $1.13\cdot 10^{-2}$  & $113436$ \\ \hline
\end{tabular}
\end{table}
\begin{figure}[H]
  \includegraphics[width=0.49\textwidth]{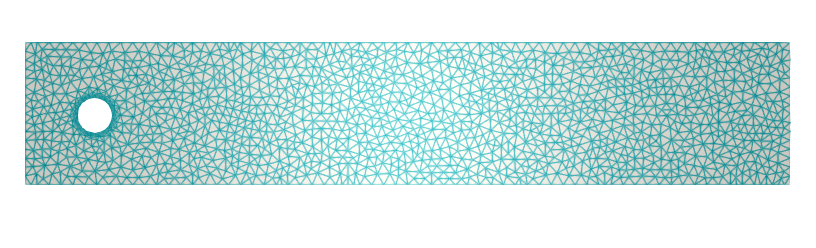}
      \includegraphics[width=0.49\textwidth]{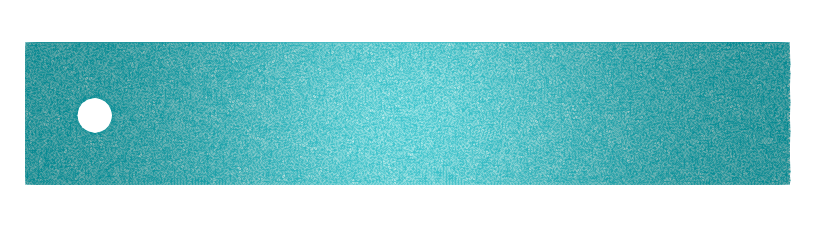} \\
    \caption{Ganeral setting. Coarse mesh (M1) and fine mesh (M2), left and right plots, respectively.}
  \label{fig:mesh}
\end{figure}

\begin{figure}[H]
  \includegraphics[width=0.49\textwidth]{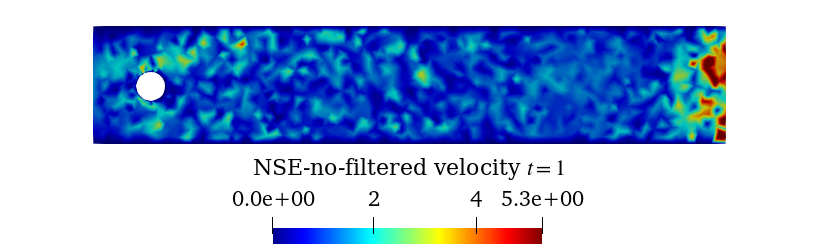}
      \includegraphics[width=0.49\textwidth]{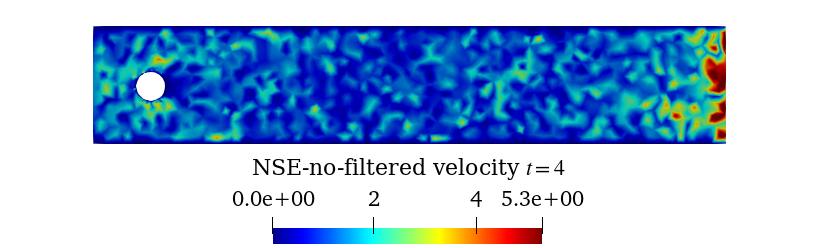} \\
    \caption{General setting. Velocity profiles for NSE-no-filtered simulation on the coarse mesh (M1) at $t=1$ and $t=4$, left and right plots, respectively.}
  \label{fig:noEF}
\end{figure}
The EF, VMS-EFFC, and VMS-EPFC algorithms are applied on a triangular mesh (left plot in Figure \ref{fig:mesh}) with
 $h_{min} = 4.46 \cdot 10^{-3}$  and $h_{max} = 4.02 \cdot 10^{-2}$ as the minimum and the maximum diameter of the elements, respectively. We employ a Taylor-Hood $\mathbb P^2 - \mathbb P^1$ FE approximation for velocity and pressure, respectively, to obtain a FE space of dimension $N_h \doteq N_h^{u} + N_h^p = 14116$. The mesh cannot accurately capture the physical flow features for high $Re$, and displays spurious numerical oscillations. This is depicted in Figure \ref{fig:noEF}, which shows the velocity profile obtained for $\nu = 10^{-4}$ (i.e., $Re=1000$) and $\Delta t = 4\cdot 10^{-4}$, at time instances $t=1$ and $t=4$, left and right plots, respectively. It is clear that stabilization strategies are needed in this context to alleviate the numerical oscillations and to yield more accurate results. \\ 
 In the next section, we present the results for EF, VMS-EFFC, and VMS-EPFC applied to the proposed setting. We recall that when a differential filter is applied, a grad-div stabilization is performed with $\gamma_D=100$ \cite{Strazzullo20223148}.
 
\subsection{Experiment 1 (flow past a cylinder): What Is The Preferable Strategy At The FOM Level?}
\begin{figure}[H]
  \includegraphics[width=0.49\textwidth]{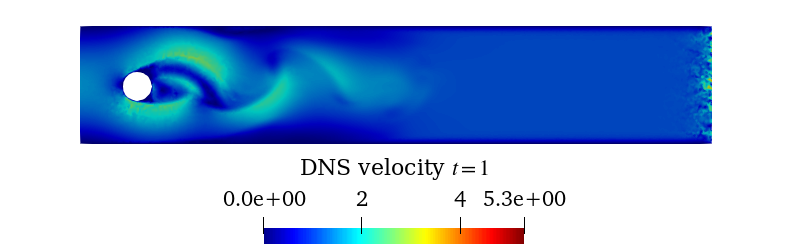}
      \includegraphics[width=0.49\textwidth]{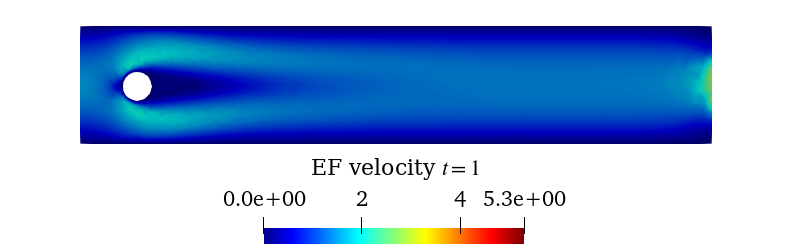} \\
    \includegraphics[width=0.49\textwidth]{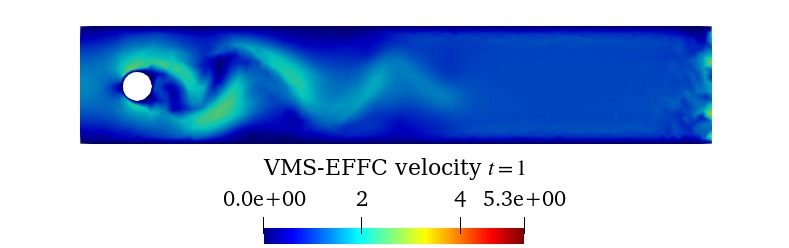}
  \includegraphics[width=0.49\textwidth]{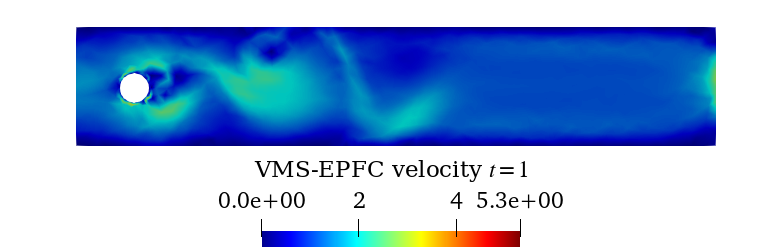}\\
  \caption{Experiment 1. Velocity profiles at $t=1$: DNS (top left), EF (top right) for $\delta=1.59 \cdot 10^{-3}$, VMS-EFFC (bottom left) for $\delta_1=\delta_2=1.59 \cdot 10^{-3}$, and VMS-EPFC (bottom right) for $\delta=1.59 \cdot 10^{-3}$.}
  \label{fig:vt1}
\end{figure}

\begin{figure}[H]
  \includegraphics[width=0.49\textwidth]{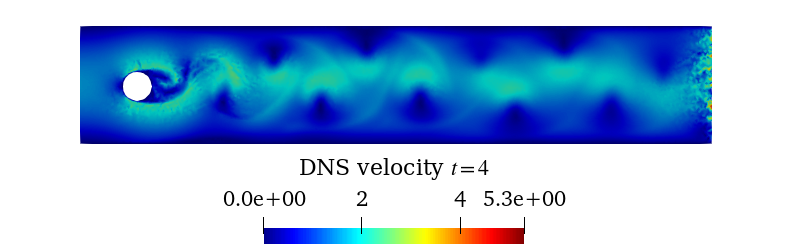}
      \includegraphics[width=0.49\textwidth]{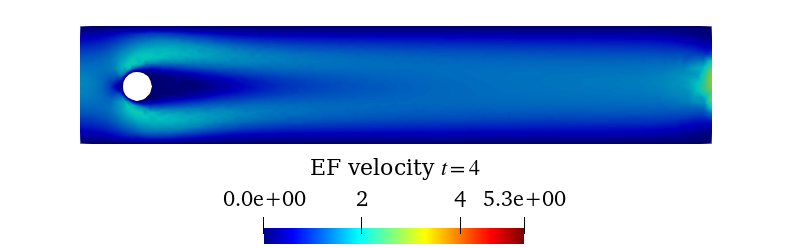} \\
    \includegraphics[width=0.49\textwidth]{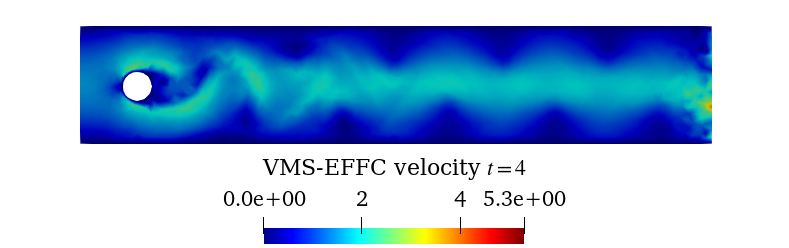}
  \includegraphics[width=0.49\textwidth]{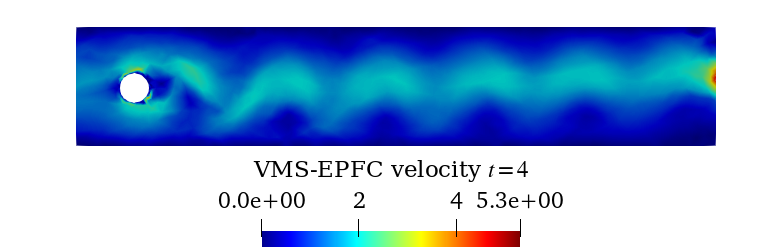}\\
 \caption{Experiment 1. Velocity profiles at $t=4$: DNS (top left), EF (top right) for $\delta=1.59 \cdot 10^{-3}$, VMS-EFFC (bottom left) for $\delta_1=\delta_2=1.59 \cdot 10^{-3}$, and VMS-EPFC (bottom right) for $\delta=1.59 \cdot 10^{-3}$.}
  \label{fig:vt4}
\end{figure}

\begin{figure}[H]
  \includegraphics[width=0.49\textwidth]{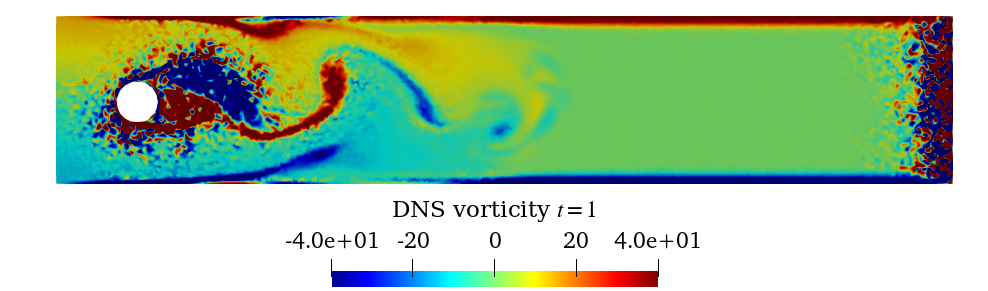}
      \includegraphics[width=0.49\textwidth]{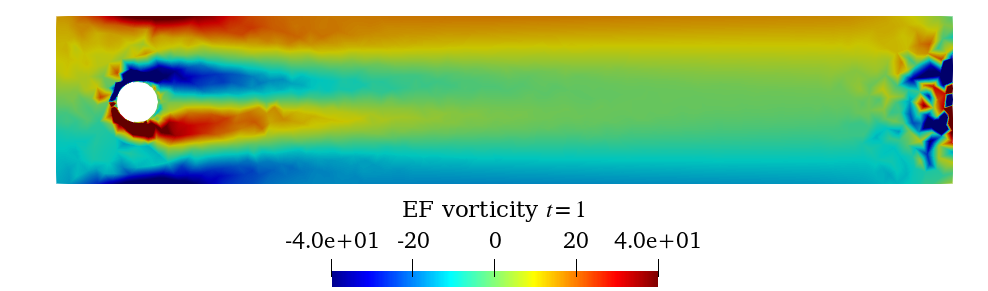} \\
    \includegraphics[width=0.49\textwidth]{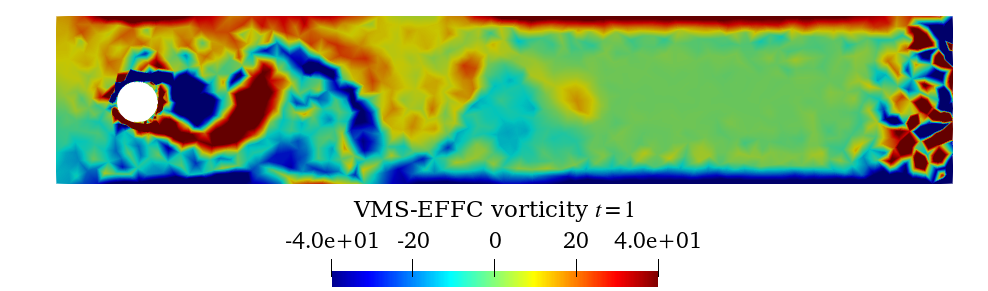}
  \includegraphics[width=0.49\textwidth]{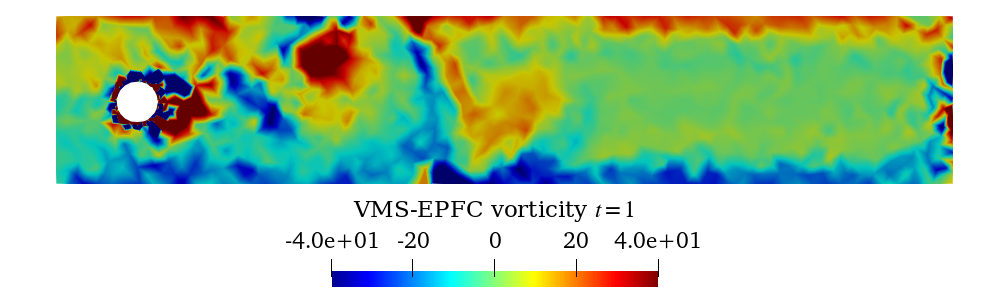}\\
  \caption{{Experiment 1. Vorticity profiles at $t=1$}\color{black}: DNS (top left), EF (top right) for $\delta=1.59 \cdot 10^{-3}$, VMS-EFFC (bottom left) for $\delta_1=\delta_2=1.59 \cdot 10^{-3}$, and VMS-EPFC (bottom right) for $\delta=1.59 \cdot 10^{-3}$.}
  \label{fig:vort1}
\end{figure}

\begin{figure}[H]
  \includegraphics[width=0.49\textwidth]{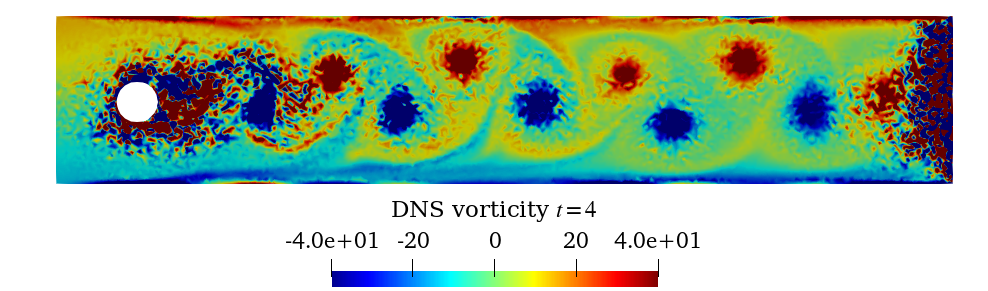}
      \includegraphics[width=0.49\textwidth]{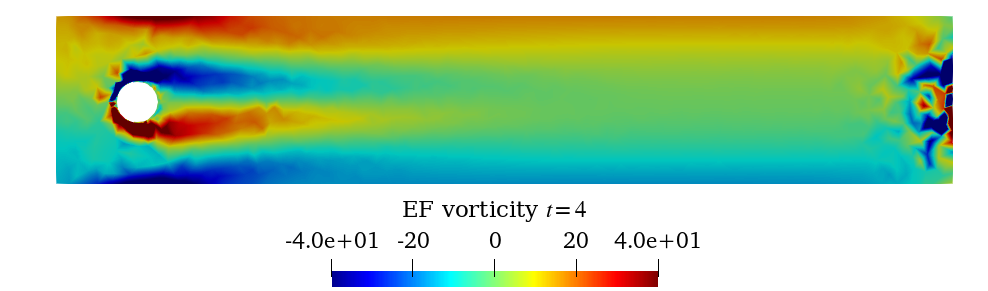} \\
    \includegraphics[width=0.49\textwidth]{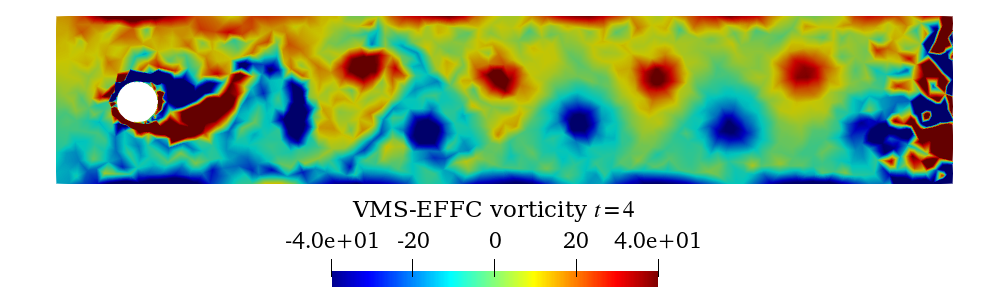}
  \includegraphics[width=0.49\textwidth]{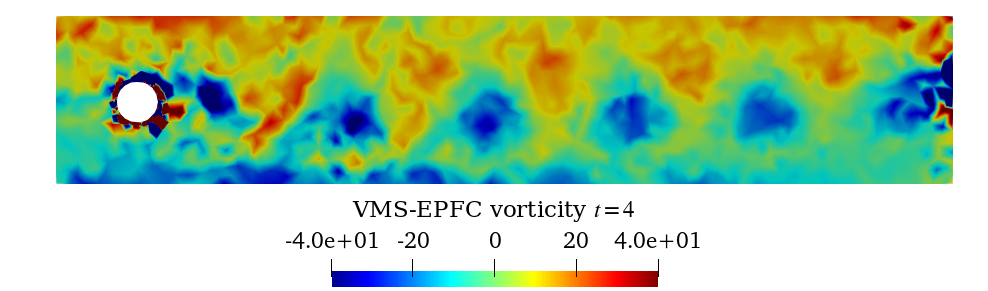}\\
  \caption{{Experiment 1. Vorticity profiles at $t=4$}\color{black}: DNS (top left), EF (top right) for $\delta=1.59 \cdot 10^{-3}$, VMS-EFFC (bottom left) for $\delta_1=\delta_2=1.59 \cdot 10^{-3}$, and VMS-EPFC (bottom right) for $\delta=1.59 \cdot 10^{-3}$.}
  \label{fig:vort4}
\end{figure}

In this section, we compare the performances of EF, VMS-EFFC, and VMS-EPFC with respect to the DNS results of NSE on a fine mesh (depicted in the right plot of Figure \ref{fig:mesh}), denoted as M2. As already specified, the filters are applied on a coarser mesh, M1. Table \ref{tab:meshval} shows the parameters of the two meshes after a Taylor-Hood $\mathbb P^2-\mathbb P^1$ discretization. The fine mesh M2 has $N_h =113436$ degrees of freedom, with $h_{min} = 4.44\cdot 10^{-3}$ and  $h_{max} = 1.13\cdot 10^{-2}$. For the sake of clarity, we recall that we are working with $\nu = 10^{-4}$, $T=4$, and $\Delta t = 4\cdot 10^{-4}$. Furthermore, we recall that every time a DF is applied, the grad-div stabilization is performed with $\gamma_D=100$. Moreover, in the postprocessing step of VMS-EPFC, the grad-div stabilization is performed with $\gamma_P = 0.01$. We refer the reader to Remark \ref{rem:gd_param} for a more detailed discussion of the grad-div stabilization effect on the VMS-based filters.\\
First, we describe the flow dynamics by an EF simulation for $\delta = 1.59 \cdot 10^{-3}$. 
The choice of the filter radius is made to showcase the benefits of using VMS-based filters in the context of overdiffusive DF. The overdiffusive action is depicted in the top right plots of Figures \ref{fig:vt1} and \ref{fig:vt4}, where the velocity fields for the EF algorithms are shown for $t=1$ and $t=4$, respectively. We compare them to the DNS velocity shown in the top left plots of Figures \ref{fig:vt1} and \ref{fig:vt4}: It is clear that the DF smooths the spurious oscillations of the NSE-no-filtered velocity of Figure \ref{fig:noEF}, but its action totally damps out the vortex shedding features of the flow. The VMS-EFFC and VMS-EPFC algorithms recover the vortex shedding. Indeed, for the overdiffusive value of $\delta = 1.59 \cdot 10^{-3}$ used in the EF simulations, the VMS-based strategies show the classical vortex shedding behavior of the test case. We plot two representative solutions at $t=1$ and $t=4$ in the bottom plots of Figures \ref{fig:vt1} and \ref{fig:vt4}: the VMS-EFFC and VMS-EPFC velocities are displayed in the left and right plots, respectively. From the plot comparison, the VMS-EFFC outperforms the other strategies, capturing the vortex shedding behavior, and avoiding unphysical behavior and overdiffusivity, which are displayed by the VMS-EPFC velocity for $t=1$ and $t=4$, respectively. {This conclusion can also be drawn by observing Figure \ref{fig:vort1} and Figure \ref{fig:vort4}, representing the vorticity for DNS, EF, VMS-EFFC, and VMS-EPFC at $t=1$ and $t=4$. The VMS-EFFC is the only strategy capable of recovering the vortex shedding behavior, which is unphysical for the VMS-EPFC approach and completely smoothed out by the EF strategy. Since the qualitative behavior of the vorticity plots and the velocity profiles is the same, in the next experiments we will only focus on the velocity profiles. }\color{black}
To quantitatively test the accuracy of the methods, we start with average quantities, well-suited to analyze regularized approaches \cite{BIL05,rebollo2014mathematical,sagaut2006large} such as the $L^2-$norm of the velocity, and the drag and lift coefficients, defined as
\begin{equation}
\label{eq:drag}
C_D(t) \doteq \frac{2}{U^2L}\int_{\partial \Omega_C}((2\nu\nabla u - p{I}) \cdot {n}_C)\cdot  {t}_C \;ds, 
\end{equation}
and
\begin{equation}
\label{eq:lift}
C_L(t) \doteq \frac{2}{U^2 L}\int_{\partial \Omega_C}((2\nu\nabla u - p{I}) \cdot {n}_C)\cdot {n}_C \;ds,
\end{equation}
for a generic velocity $u$ and pressure $p$, where ${n}_C$ and ${t}_C$ are the normal and tangential unit vectors to the cylinder boundary, respectively.
We also use the relative error for velocity fields
{at} each time instance:
$$
E_u(t^n) = \frac{\norm{\widetilde{u}^{n} - u^n}_{L^2(\Omega)}}{\norm{\widetilde{u}^{n}}_{L^2(\Omega)}},
$$ where $u^n$ is the EF, VMS-EFFC, and VMS-EPFC velocity at time $t_n$, and $\widetilde{u}^{n}$ is the projection of the DNS results on the coarser mesh.

\begin{figure}[H]
  \centering
  \includegraphics[width=0.45\textwidth]{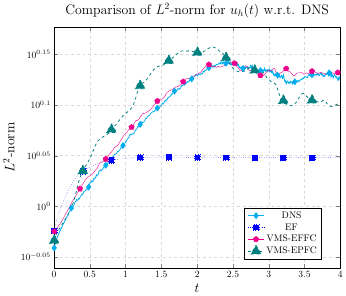}
  \includegraphics[width=0.45\textwidth]{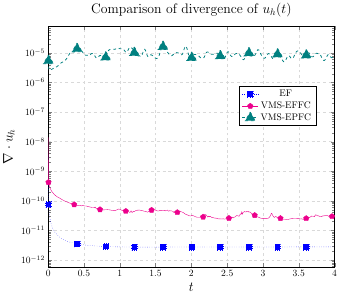}
  \caption{Experiment 1. \emph{Left}: Velocity $L^2$-norms in time for EF with $\delta=1.59 \cdot 10^{-3}$, VMS-EFFC for $\delta_1=\delta_2=1.59 \cdot 10^{-3}$, and VMS-EPFC for $\delta=1.59 \cdot 10^{-3}$ compared to the DNS solution. \emph{Right}: divergence evolution for EF with $\delta=1.59 \cdot 10^{-3}$, VMS-EFFC with $\delta_1=\delta_2=1.59 \cdot 10^{-3}$, and VMS-EPFC with $\delta=1.59 \cdot 10^{-3}$.}
  \label{fig:norms}
\end{figure}

Figure \ref{fig:norms} shows the squared $L^2$-norms of the EF, VMS-EFFC, and VMS-EPFC velocities compared to the squared $L^2$-norm of the DNS velocity projected on the coarser mesh, and the evolution of the divergence of the three algorithms, in the left and right plots, respectively. It is evident that the EF velocity does not recover the DNS energy norm. In contrast, VMS-EFFC and VMS-EPFC yield much more accurate results. Furthermore, VMS-EFFC is significantly more accurate than VMS-EPFC. 

\begin{figure}[H]
  \centering
  \includegraphics[width=0.486\textwidth]{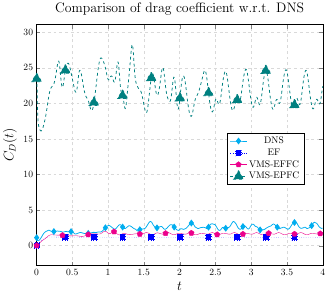}
  \includegraphics[width=0.49\textwidth]{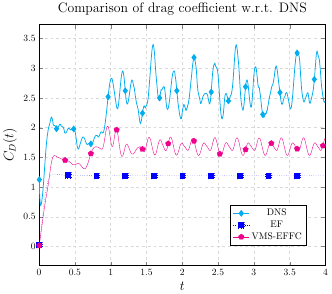}
  \caption{Experiment 1. $C_D(t)$: DNS, EF for $\delta=1.59 \cdot 10^{-3}$, VMS-EFFC for $\delta_1=\delta_2=1.59 \cdot 10^{-3}$, and VMS-EPFC for $\delta=1.59 \cdot 10^{-3}$, left plot. In the right plot, we show the same quantities but we exclude the VMS-EPFC from the results.}
  \label{fig:drag}
\end{figure}
\begin{figure}[H]
  \centering
  \includegraphics[width=0.488\textwidth]{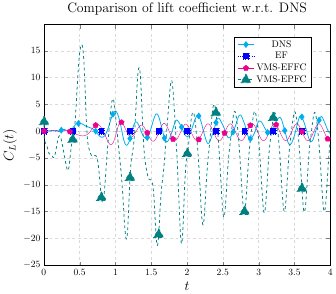}
  \includegraphics[width=0.49\textwidth]{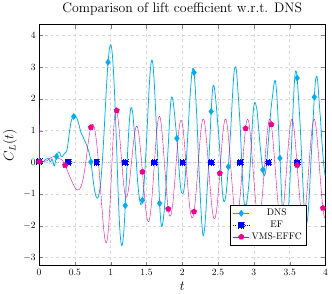}
  \caption{Experiment 1. $C_L(t)$: DNS, EF for $\delta=1.59 \cdot 10^{-3}$, VMS-EFFC for $\delta_1=\delta_2=1.59 \cdot 10^{-3}$, and VMS-EPFC for $\delta=1.59 \cdot 10^{-3}$, left plot. In the right plot, we show the same quantities but we exclude the VMS-EPFC from the results.}
  \label{fig:lift}
\end{figure}

\begin{figure}[H]
  \centering
  \includegraphics[width=0.45\textwidth]{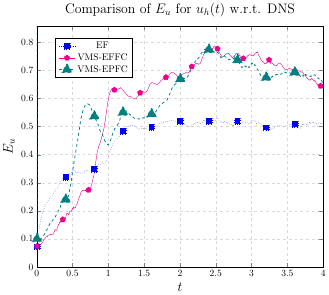}
  \caption{Experiment 1. Velocity relative error in time for EF with $\delta=1.59 \cdot 10^{-3}$, VMS-EFFC for $\delta_1=\delta_2=1.59 \cdot 10^{-3}$, and VMS-EPFC for $\delta=1.59 \cdot 10^{-3}$ compared to the DNS solution.}
  \label{fig:errs}
\end{figure}

\begin{figure}[H]
  \centering
  \includegraphics[width=0.32\textwidth]{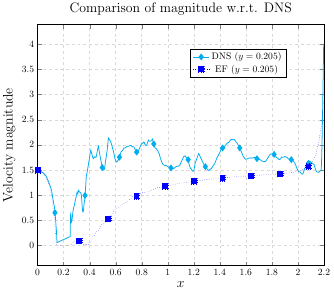}
  \includegraphics[width=0.32\textwidth]{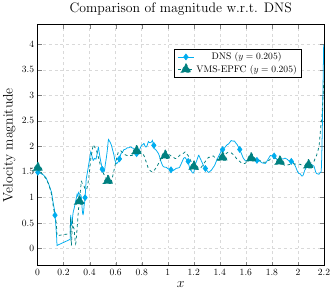}
  \includegraphics[width=0.32\textwidth]{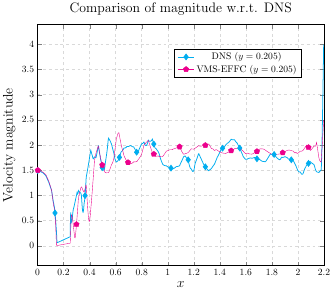}\\

    \includegraphics[width=0.32\textwidth]{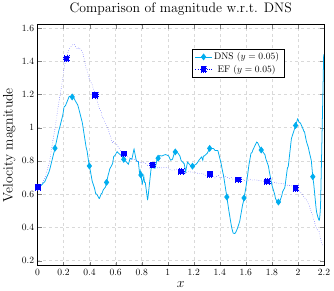}
  \includegraphics[width=0.32\textwidth]{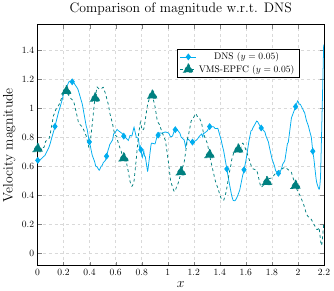}
  \includegraphics[width=0.32\textwidth]{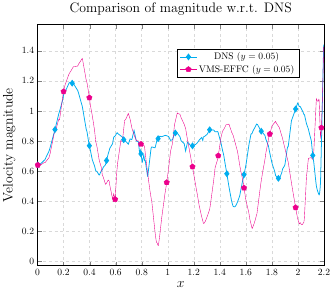}\\

      \includegraphics[width=0.32\textwidth]{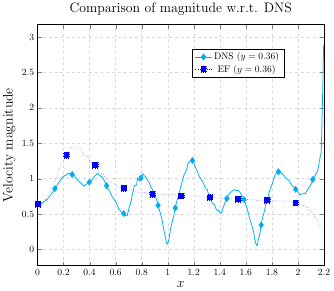}
  \includegraphics[width=0.32\textwidth]{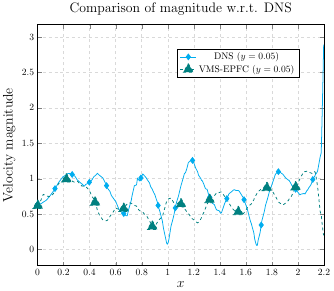}
  \includegraphics[width=0.32\textwidth]{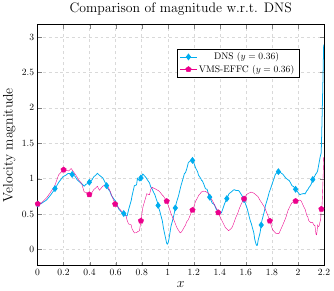}
  
  \caption{Experiment 1. Comparison of the velocity for fixed y values and varying x values.  DNS results are shown for comparison purposes. Top row: EF for $\delta=1.59 \cdot 10^{-3}$ (left), VMS-EPFC for $\delta_1=\delta_2=1.59 \cdot 10^{-3}$ (center), and VMS-EFFC for $\delta=1.59 \cdot 10^{-3}$ (right) for $y=0.205$. Center row: EF for $\delta=1.59 \cdot 10^{-3}$ (left), VMS-EPFC for $\delta_1=\delta_2=1.59 \cdot 10^{-3}$ (center), and VMS-EFFC for $\delta=1.59 \cdot 10^{-3}$ (right) for $y=0.05$.
  Bottom row: EF for $\delta=1.59 \cdot 10^{-3}$ (left), VMS-EPFC for $\delta_1=\delta_2=1.59 \cdot 10^{-3}$ (center), and VMS-EFFC for $\delta=1.59 \cdot 10^{-3}$ (right) for $y=0.36$.}
  \label{fig:overline}
\end{figure}

In Figures \ref{fig:drag} and \ref{fig:lift}, we show the drag and lift coefficients for the EF, VMS-EFFC, and VMS-EPFC in the left plots. The coefficients are compared to the DNS coefficients. The plots for $C_D(t)$ and $C_L(t)$ show that EF is overdiffusive: the magnitude of its coefficients is generally much lower than the magnitude of the DNS coefficients. Moreover, VMS-EPFC fails to recover the coefficients, while VMS-EFFC performs better in terms of the magnitude of the coefficients with respect to the DNS simulation. {However, in Figure \ref{fig:lift}, we can observe a phase shift for the VMS-EFFC. This will affect the accuracy of the solution in terms of relative error; we postpone more detailed considerations to the end of the section.}
\\ In the right plots of Figures \ref{fig:drag} and \ref{fig:lift}, we exclude VMS-EPFC from the representation to better compare the performances of the other strategies. The poor results of the VMS-EPFC algorithm for the drag and lift coefficients are due to the pressure field representation, which is presented in Appendix \ref{app:pressure}: the VMS-EPFC fails in reconstructing the pressure, while VMS-EFFC yields more accurate results. The divergence value of the solution represents another drawback of the VMS-EPFC. In the right plot of Figure \ref{fig:norms}, we show the evolution of the divergence of EF, VMS-EFFC, and VMS-EPFC. We observe that the first two strategies yield divergence values on the order of $10^{-11}$ and $10^{-10}$, while the VMS-EPFC reaches divergence values higher than $10^{-5}$.
For the sake of completeness, we show the \textit{pointwise} relative error in time for all three approaches in Figure \ref{fig:errs}. From the plot, we notice that all the strategies show large $E_u$ values, and EF looks to be the best approach. This behavior is somewhat unexpected since the VMS-EFFC outperforms the other methods in recovering the $L^2-$norm. However, this unexpected behavior can be explained in Figure \ref{fig:overline} where we show the velocity magnitude of EF, VMS-EPFC, and VMS-EFFC for $t=4$, for all $x \in [0,2,2]$ and for a fixed $y$ values, i.e., the magnitude for \emph{horizontal lines} in the spatial domain at the final time. In the top row, we can observe how VMS-EFFC outperforms the other two strategies, closely following the DNS behavior for $y=0.205$, i.e., in the middle of the spatial domain.

As expected the EF strategy is overdiffusive, while VMS-EFFC and VMS-EPFC approaches yield more accurate results. We draw the same conclusions observing the bottom row in Figure \ref{fig:overline}, where we fix $y=0.36$. Namely, at the top of the domain, VMS -EFFC and VMS-EPFC strategies give a better representation of the DNS velocity compared to the standard EF.
For the middle row in Figure \ref{fig:overline}, which corresponds to the bottom of the spatial domain for $y=0.05$, none of the strategies can recover the DNS velocity. The VMS-EFFC field presents large oscillations in this case, and this yields large values of the relative error. However, for $y=0.205$ and $y=0.36$ is the most accurate strategy.\\
The velocity in the bottom part {of the domain} is more complicated to represent. {This is already visible} in  Figure \ref{fig:vt4}, {where the DNS behavior is more complex at the bottom of the domain}\footnote{We are using unstructured meshes, and this may lead to asymmetries in the configuration.}, with several recirculations areas arising. The VMS-EFFC and VMS-EPFC strategies are smoothing these vortices and these features are lost in the simulations. Nevertheless, the VMS-EFFC and the VMS-EPFC results present a phase shift, {as already evidenced when discussing} Figure \ref{fig:lift}. {The phase shift is also evident in the vorticity plots in Figures \ref{fig:vort1} and \ref{fig:vort4}}\color{black}. This causes an increase in the pointwise error.
From these considerations, we conclude that the relative error is not an appropriate metric for testing the accuracy of the VMS-based filters for this specific test case, while average quantities such as the $L^2$-norm are more appropriate to analyze the performances of regularized models. Of course, this is well known in the numerical simulation of turbulent flows, where averaged quantities are used to evaluate the performance of the computational models. 

In conclusion, due to the poor EF and VMS-EPFC performances, we choose VMS-EFFC as the FOM model for our analysis of VMS-based filters at the ROM level.\\
We also remark that the results are sensitive to the values of $\delta$, $\delta_1$, and $\delta_2$. We present additional results on this topic in Appendix \ref{app:diff}. 
\\
We stress that the results presented in this section suggest that VMS-based strategies can be helpful when dealing with an overdiffusive action of the filter radius. We also note that, for some values of $\delta$, the EF strategy can be as competitive as the VMS-based algorithm. However, finding the \emph{optimal} setting for the various experiments can be challenging and is beyond the scope of the manuscript. The main goal of this section was to provide some preliminary results to guide us in choosing the FOM model.  We also note that, for the sake of brevity, we showed only the velocity field results. Additional results on pressure reconstruction can be found in Appendix \ref{app:pressure}. \\
To summarize, the VMS-EFFC simulations are the closest to the DNS simulations in terms of overall qualitative behavior and $L^2$-norm of the velocity. We note, however, that none of the approaches we tested can accurately describe the velocity field in terms of drag and lift coefficients, and relative errors.

\begin{remark}[More on the relative error]
There are several ways to decrease the VMS-EFFC and VMS-EPFC errors with respect to the DNS simulations, such as decreasing $\Delta t$, increasing the filter radius, and decreasing the Reynolds number. In our numerical examples, these approaches yielded relative errors that were comparable to the EF relative error, however, the average quantities suffered and the qualitative behavior was worse than the one illustrated in Figures \ref{fig:vt1} and \ref{fig:vt4}. The proposed test case is a trade-off between accuracy and computational costs effort for $Re=1000$ and represents a ground setting for the ROMs results of Section \ref{sec:exp2}.
\end{remark}

\begin{remark}[Choosing the grad-div stabilization parameter $\gamma_P$]
\label{rem:gd_param} 
For the DF, the choice of the grad-div stabilization parameters has already been investigated (see, e.g., \cite{Strazzullo20223148} and the references therein). We emphasize, however, that the grad-div parameter in step (II) of the VMS-EPFC strategy has not been investigated. We take a first step in this direction and perform a numerical investigation of the grad-div parameter in step (II) of VMS-EPFC, which yields an optimal value $\gamma_P = 0.01$. As expected, higher values of $\gamma_P$ give overdiffusive  and unphysical results, while smaller values of $\gamma_P$ result in unacceptable values for the divergence of the velocity field. For the sake of brevity, we do not include these results.
\end{remark}

\begin{remark}[Relaxation]
In the literature, the evolve-filter-relax (EFR) strategy has been used to tackle the overdiffusive action that may arise in the EF approach. The EFR strategy applies Algorithm \ref{alg:EF}, but the new velocity is given by a convex combination of the evolved velocity field, $w^{n+1}$, and the filtered velocity field, $\overline{w}^{n+1}$. Namely, for an appropriately chosen relaxation parameter $0 \leq \chi \leq 1$, the new velocity is
$$
u^{n+1} = (1 - \chi)w^{n+1} + \chi \overline{w}^{n+1}.
$$
We note that, in our numerical setting, the relaxation is not a practical option since it adds spurious oscillations to the final results. In Figure \ref{fig:EFR}, we show the velocity field obtained using EFR with $\delta = 1.59 \cdot 10^{-3}$ and $\chi = \Delta t$ (the choice of $\chi$ is discussed in \cite{bertagna2016deconvolution}). 
\begin{figure}[H]
  \includegraphics[width=0.49\textwidth]{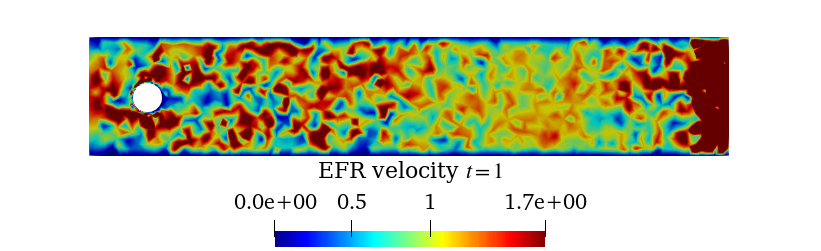}
      \includegraphics[width=0.49\textwidth]{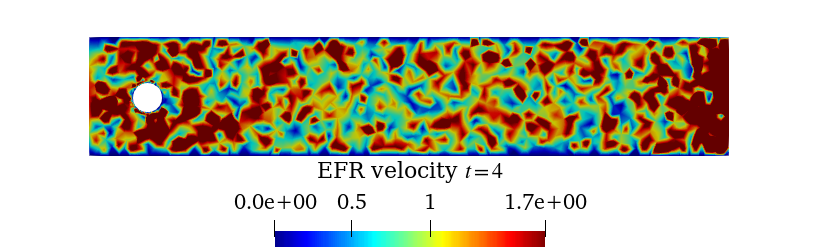} \\
    \caption{Experiment 1. Velocity profiles for EFR simulation on the coarse mesh (M1) for $t=1$ and $t=4$, left and right plots, respectively.}
  \label{fig:EFR}
\end{figure}
\end{remark}

{
\begin{remark}[Robustness of the VMS-EFFC approch]
To test the robustness of the VMS-EFFC approach, i.e., the best approach for the flow past a cylinder case, we propose another experiment: the lid-driven cavity flow at $Re=7500$. We consider the domain $\Omega=(0,1)^2$. The problem features the following Dirichlet boundary conditions\footnote{To guarantee the uniqueness of $p$, we fix its value to zero in the origin of the system coordinates.}: 
\begin{equation}
\label{eq:inlet}
u_D = 
\begin{cases}
0 \hspace{3.8cm} \text{ in }\Gamma_{\text{w}}, \\
u_{\text{in}} = \begin{cases}
\displaystyle \left  ( \frac{10x}{0.06}, 0 \right) & \text{ for $x \leq 0.6$ in } \Gamma_{\text{in}}, \\
\displaystyle \left  ( \frac{10(1-x)}{0.06}, 0 \right) & \text{ for $x \geq 0.96$ in } \Gamma_{\text{in}}, \\
\displaystyle 10 & \text{ elsewhere in } \Gamma_{\text{in}}, 

\end{cases}
\end{cases}
\end{equation}
where $\Gamma_{\text{in}} = [0,1] \times \{1\}$ is the top boundary of the unit square domain, and $\Gamma_\text{w} = \partial \Omega \setminus \Gamma_{\text{in}}$.
The initial condition is $u_0 = (0,0)$. 
We compare the performances of EF and VMS-EFFC with respect to the DNS results of the NSE on a fine mesh denoted as M4. As usual, the filters are applied on a coarser mesh, M3. Table \ref{tab:meshval_2} collects the parameters of the two meshes after a Taylor-Hood $\mathbb P^2-\mathbb P^1$ discretization. The fine mesh M4 has $N_h =140099$ degrees of freedom, with $h_{min} = 6.58 \cdot 10^{-3}$ and  $h_{max} = 1.26\cdot 10^{-2}$.
The coarse mesh M3 has $N_h =14635$ degrees of freedom, with $h_{min} = 2.05 \cdot 10^{-2}$ and  $h_{max} = 3.97\cdot 10^{-2}$.  We notice that the mesh parameters of M3 and M4 are comparable to M1 and M2 of Experiment 1, respectively. Specifically, M1 presents a smaller $h_{min}$, but this is only due to a refinement in a neighborhood of the cylinder. 
We are using $\nu = 2.5\cdot 10^{-3}$, $T=2$, and $\Delta t = 4\cdot 10^{-4}$. Also for this test case, every time a DF is applied, the grad-div stabilization is performed with $\gamma_D=100$ and $\delta = \delta_1 = \delta_2 = 1.59 \cdot 10^{-3}$.
We plot two representative solutions at $t=1$ and $t=2$ in Figures \ref{fig:vt1_c} and \ref{fig:vt4_c} for DNS, VMS-EFFC, and EF. From the plots comparison, we observe that the VMS-EFFC outperforms EF, capturing more accurately the vortex behavior, and avoiding the overdiffusivity that characterizes the EF velocity for both time instances.
To quantitatively test the accuracy of the methods, we analyze the $L^2$-norm of the velocity and the relative error with respect to the DNS simulation in Figure \ref{fig:err_c}. In the left plot, we compare the squared $L^2$-norms of the EF and the VMS-EFFC with the squared $L^2$-norm of the DNS velocity projected on the coarser mesh. The EF velocity does not recover the DNS energy norm, while the VMS-EFFC yields much more accurate results: the DNS and the VMS-EFFC norms almost coincide. The right plot shows that, for this test case, the VMS-EFFC strategy also outperforms EF in terms of the velocity relative error. Nevertheless, the VMS-EFFC yields large values of the relative error. This is due to a slight shift of the vortex to the bottom of the cavity at the beginning of the simulation. This phenomenon is observable in Figure \ref{fig:vt1_c}, where the solution is depicted for $t=1$. As already noticed, the relative error is not an appropriate \emph{general} metric for testing the accuracy of regularized models, while average quantities such as the $L^2$-norm are a robust indicator to assess the performances of regularized models.
Furthermore, we remark that:
\begin{itemize}
    \item[$\circ$] If no regularization is applied, the Newton solver for NSE over the coarse mesh does not converge. Thus, regularization is needed for this test case.
    \item[$\circ$] To guarantee a fair comparison with respect to Experiment 1, for the lid-driven cavity test, we keep the same parameters for DFs and grad-div stabilization used for the flow past a cylinder test.
    \item[$\circ$] Grad-div stabilization ensures divergence values of the order of $10^{-11}$ both for EF and VMS-EFFC for all the time evolution. For the sake of brevity, we do not show these plots.
    \item[$\circ$] We also analyzed the VMS-EPCF strategy for this test case. However, VSM-EFFC proved to be the best approach in terms of all the metrics we investigated. Thus, for the sake of brevity, we only reported the VMS-EFFC performances compared to EF.  
\end{itemize}
\begin{table}[H]
\caption{General setting (lid-driven cavity). Comparison of the coarse mesh (M3), and the fine mesh (M4) with respect to $h_{min}$, $h_{max}$, and number of 
{degrees of freedom}.}
\vspace{3mm}
\label{tab:meshval_2} \centering
\begin{tabular}{|c|c|c|c|}
\hline
mesh & $h_{min}$            & $h_{max}$            &  $N_h$     \\ \hline
M3   & $2.05\cdot 10^{-2}$  & $3.97\cdot 10^{-2}$ & $14635$  \\ \hline
M4   & $6.58 \cdot 10^{-3}$ & $1.26\cdot 10^{-2}$  & $140099$ \\ \hline
\end{tabular}
\end{table}

\begin{figure}[H]
  \includegraphics[width=0.32\textwidth]{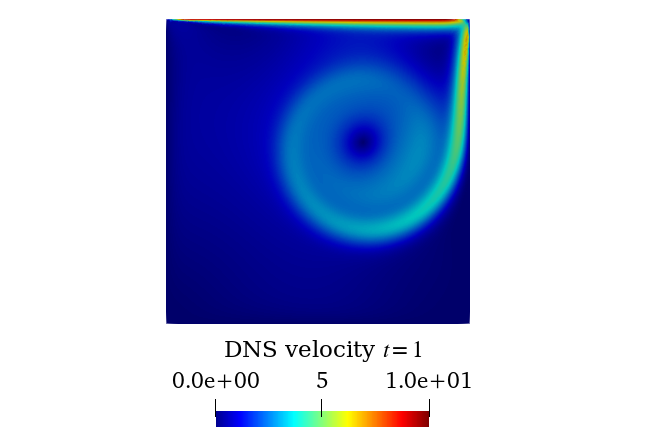}
      \includegraphics[width=0.32\textwidth]{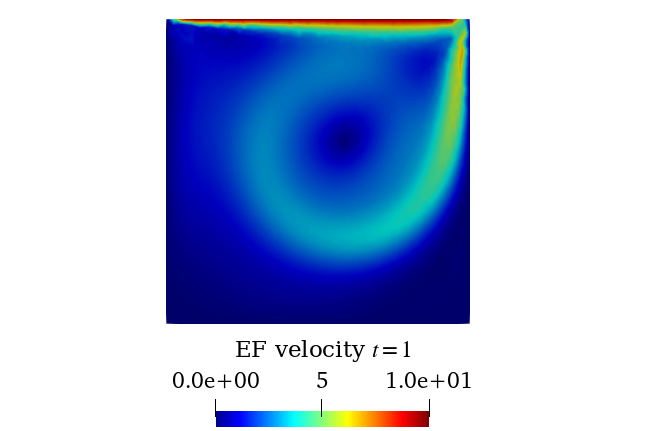}
    \includegraphics[width=0.32\textwidth]{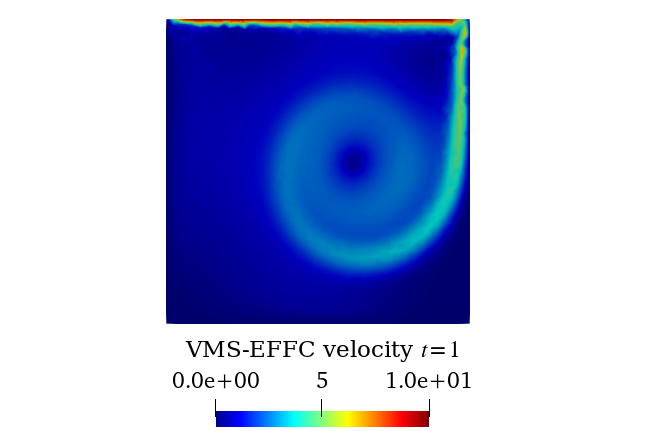}
  \caption{Experiment 1 (lid-driven cavity). Velocity profiles at $t=1$: DNS (left), EF (center) for $\delta=1.59 \cdot 10^{-3}$, and VMS-EFFC (right) for $\delta_1=\delta_2=1.59 \cdot 10^{-3}$.}
  \label{fig:vt1_c}
\end{figure}

\begin{figure}[H]
   \includegraphics[width=0.32\textwidth]{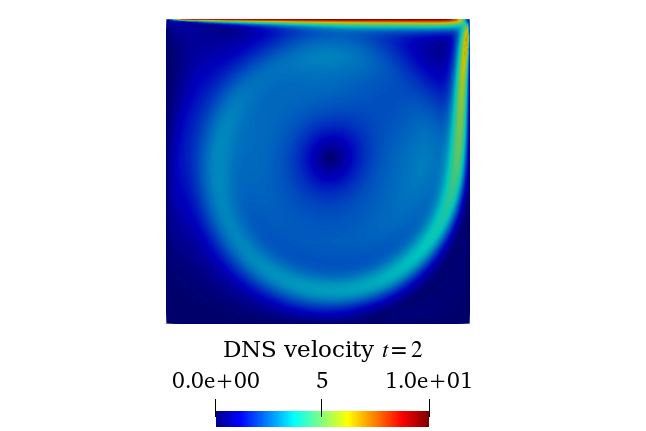}
      \includegraphics[width=0.32\textwidth]{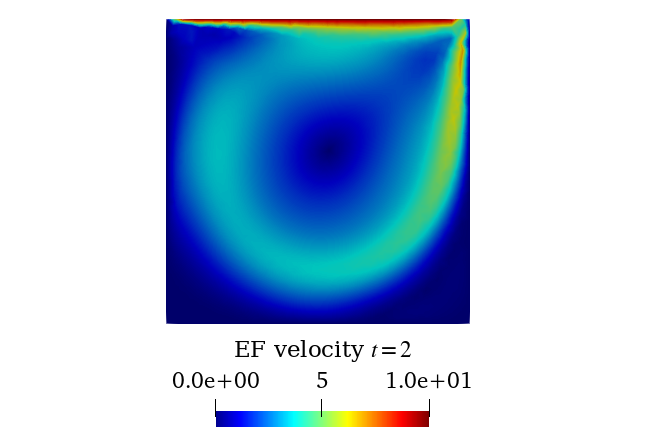}
    \includegraphics[width=0.32\textwidth]{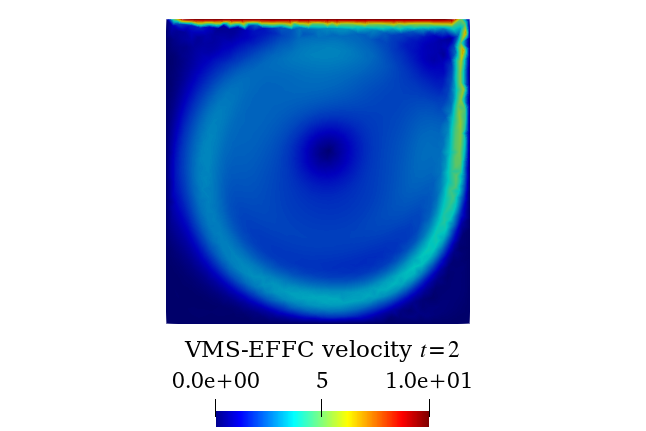}
 \caption{Experiment 1 (lid-driven cavity). Velocity profiles at $t=2$: DNS (left), EF (center) for $\delta=1.59 \cdot 10^{-3}$, and VMS-EFFC (right) for $\delta_1=\delta_2=1.59 \cdot 10^{-3}$.}
  \label{fig:vt4_c}
\end{figure}

\begin{figure}[H]
  \centering
  \includegraphics[width=0.488\textwidth]{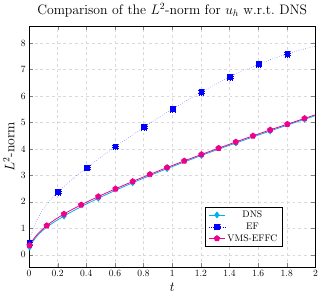}
  \includegraphics[width=0.49\textwidth]{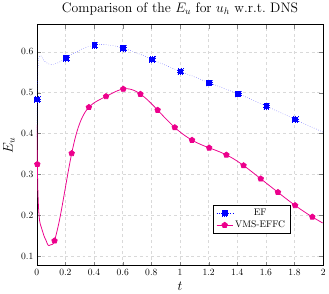}
  \caption{Experiment 1 (lid-driven cavity). \emph{Left}: Velocity $L^2$-norms in time for DNS, EF for $\delta=1.59 \cdot 10^{-3}$, and VMS-EFFC for $\delta_1=\delta_2=1.59 \cdot 10^{-3}$. \emph{Right}: Velocity relative error in time for DNS, EF for $\delta=1.59 \cdot 10^{-3}$, and VMS-EFFC for $\delta_1=\delta_2=1.59 \cdot 10^{-3}$.}
  \label{fig:err_c}
\end{figure}
\end{remark}
}\color{black}

\section{ROM Results}
\label{sec:ROM_results}
In this section, we seek the best strategy to be used at the ROM level between G-ROM, EF-ROM, VMS-EFFC-ROM and VMS-EPFC-ROM. For the sake of presentation, we focus on the results for the velocity. Additional results concerning pressure are presented in Appendix \ref{app:pressureROM}. 
\subsection{Experiment 2: What Is The Preferable Strategy At The ROM Level?}
\label{sec:exp2}
\begin{figure}[H]
\centering
  \includegraphics[width=0.49\textwidth]{Images/FOM/EFFC1.png}\\
      \includegraphics[width=0.49\textwidth]{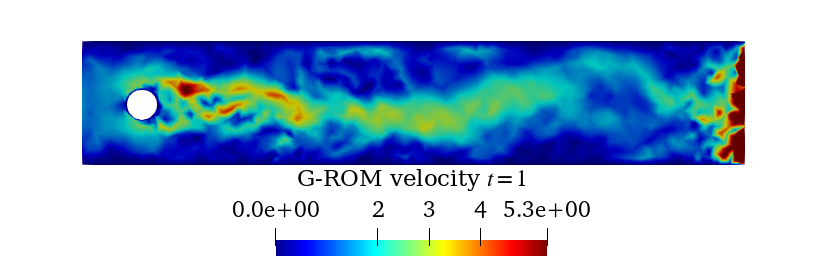}
    \includegraphics[width=0.49\textwidth]{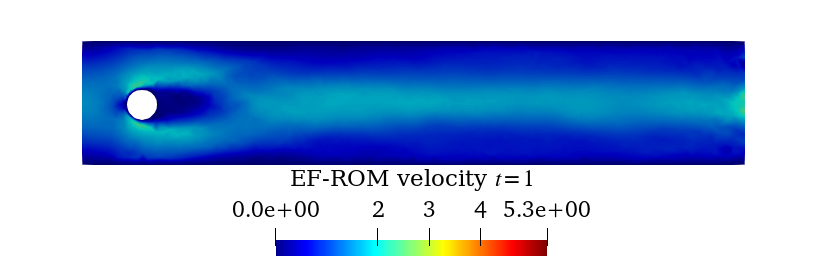}
  \includegraphics[width=0.49\textwidth]{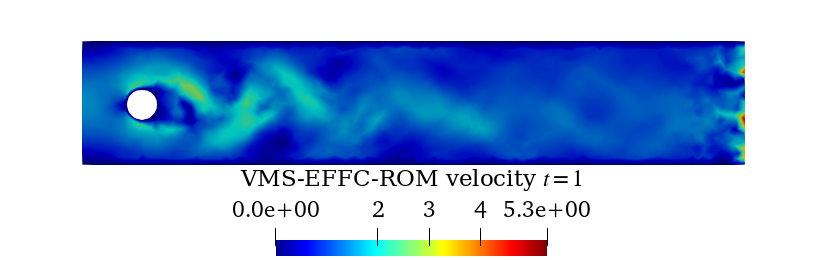}   \includegraphics[width=0.49\textwidth]{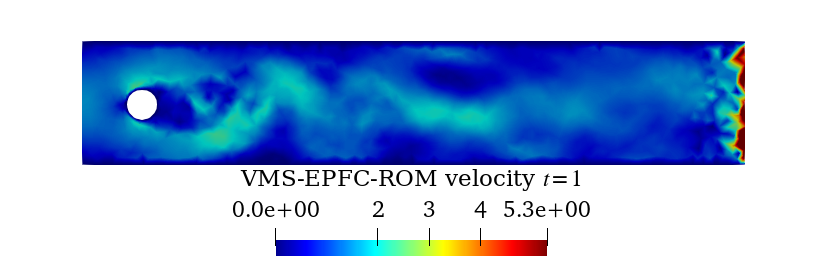}\\
  \caption{Experiment 2. Velocity profiles at $t=1$: FOM (VMS-EFFC) for $\delta=1.59 \cdot 10^{-3}$, G-ROM, EF-ROM for $\delta=1.59 \cdot 10^{-3}$, VMS-EFFC-ROM for $\delta_1=\delta_2=1.59 \cdot 10^{-3}$, VMS-EPFC-ROM for $\delta=1.59 \cdot 10^{-3}$ and $\overline{r}_{u}=77$, top, middle left, middle right, bottom left, and bottom right plots, respectively.}
  \label{fig:romvt1}
\end{figure}

\begin{figure}[H]
\centering
  \includegraphics[width=0.49\textwidth]{Images/FOM/EFFC4.png}\\
      \includegraphics[width=0.49\textwidth]{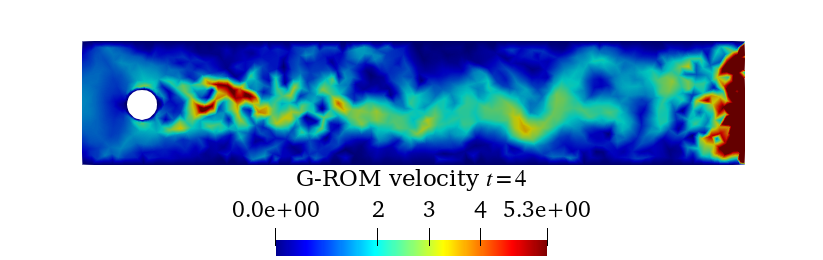}
    \includegraphics[width=0.49\textwidth]{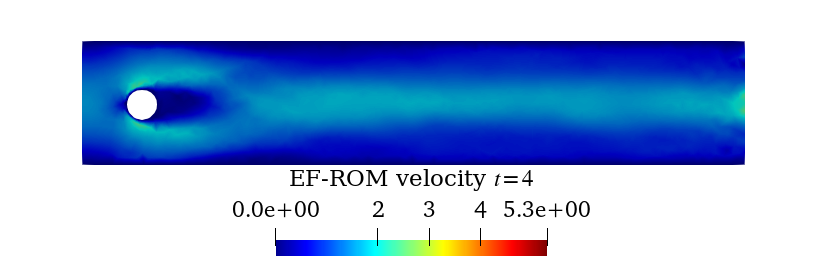}
  \includegraphics[width=0.49\textwidth]{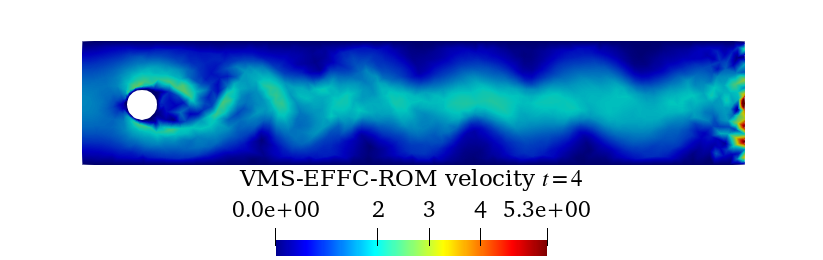}   \includegraphics[width=0.49\textwidth]{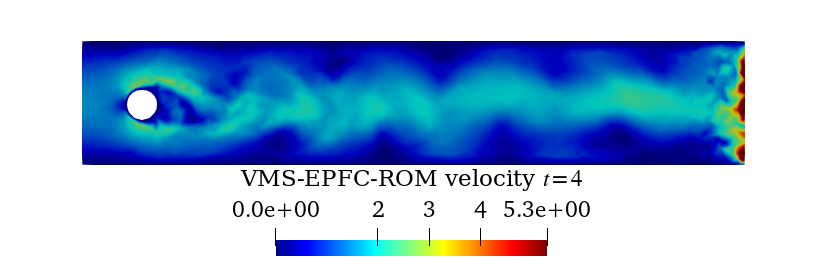}\\
  \caption{Experiment 2. Velocity profiles at $t=4$: FOM (VMS-EFFC) for $\delta=1.59 \cdot 10^{-3}$, G-ROM, EF-ROM for $\delta=1.59 \cdot 10^{-3}$, VMS-EFFC-ROM for $\delta_1=\delta_2=1.59 \cdot 10^{-3}$, VMS-EPFC-ROM for $\delta=1.59 \cdot 10^{-3}$ and $\overline{r}_{u}=77$, top, middle left, middle right, bottom left, and bottom right plots, respectively.}
  \label{fig:romvt4}
\end{figure}

\begin{figure}[H]
  \centering
  \includegraphics[width=0.45\textwidth]{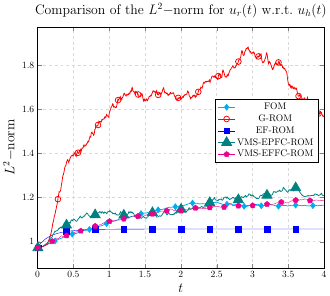}
  \includegraphics[width=0.45\textwidth]{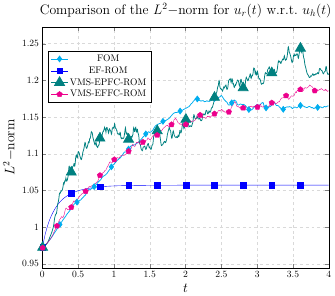}
  \caption{Experiment 2. \emph{Left}: Evolution of the squared $L^2$-norms of the velocity for G-ROM, EF-ROM for $\delta=1.59 \cdot 10^{-3}$, VMS-EFFC for $\delta_1=\delta_2=1.59 \cdot 10^{-3}$, and VMS-EPFC for $\delta=1.59 \cdot 10^{-3}$ and $\overline{r}_{u}=77$ compared to the FOM solution (VMS-EFFC) for $\delta=1.59 \cdot 10^{-3}$. The left plot shows G-ROM, which is omitted in the right plot for clarity.}
  \label{fig:romnorms}
\end{figure}

\begin{figure}[H]
  \centering
  \includegraphics[width=0.479\textwidth]{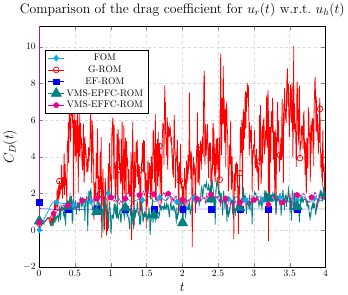}
  \includegraphics[width=0.49\textwidth]{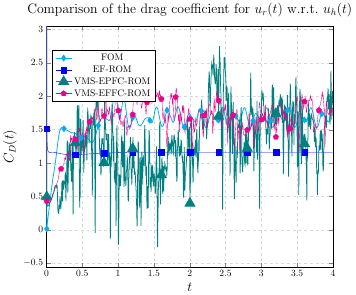}
  \caption{Experiment 2. $C_D(t)$: FOM (VMS-EFFC), for $\delta=1.59 \cdot 10^{-3}$, EF-ROM, VMS-EFFC-ROM for $\delta_1=\delta_2=1.59 \cdot 10^{-3}$, and VMS-EPFC-ROM for $\delta=1.59 \cdot 10^{-3}$ and $\overline{r}_u=77$. The left plot shows G-ROM, which is omitted in the right plot for clarity.}
  \label{fig:romdrag}
\end{figure}

\begin{figure}[H]
  \centering
  \includegraphics[width=0.43\textwidth]{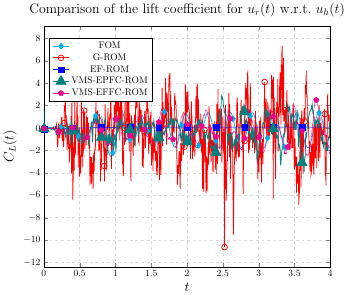}
  \includegraphics[width=0.54\textwidth]{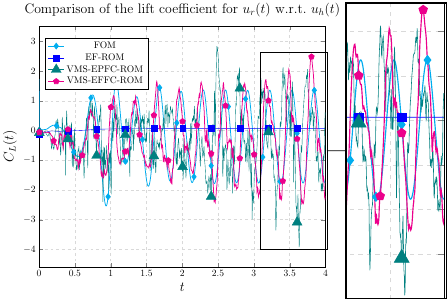}
  \caption{Experiment 2. $C_L(t)$: FOM (VMS-EFFC), for $\delta=1.59 \cdot 10^{-3}$, EF-ROM, VMS-EFFC-ROM for $\delta_1=\delta_2=1.59 \cdot 10^{-3}$, and VMS-EPFC-ROM for $\delta=1.59 \cdot 10^{-3}$ and $\overline{r}_u=77$. The left plot shows G-ROM, which is omitted in the right plot for clarity. 
  }
  \label{fig:romdlift}
\end{figure}

In this test case, we fix a FOM solution (as already specified in Section \ref{sec:FOM_results}, we use VMS-EFFC at the FOM level), and we analyzed G-ROM, EF-ROM, VMS-EFFC-ROM, and VMS-EPFC-ROM in terms of accuracy with respect to the FOM solution. The numerical setting is analogous to the one described in Section \ref{sec:FOM_results}. We perform a POD over the time interval considering $N_u=N_p=1000$ equispaced snapshots. We retain $r_u = 140$ modes for the velocity. This value preserves $99.5\%$ of the velocity energy in the system. For the pressure, we consider $r_p = 15$, which retains 99.9\% of the energy of the pressure field. We remark that 99.9\% of the velocity energy is retained for $r_u=185$, beyond the maximum number of basis functions to be considered, i.e., $N_{max}=150$. Moreover, we consider $r_s=r_p$, as suggested in \cite{ballarin2015supremizer}. Thus, we obtain $r_{us} = 155$.
In order to apply the VMS-EPFC-ROM, we set $\overline{r}_{u} = \lfloor r_u/2 \rfloor = 77$, where $\lfloor \cdot \rfloor$ represent the flooring function. We postpone the rationale behind this choice to Remark \ref{rem:role_proj}.
At each time step, to test the accuracy of the various methods, we compare the FOM and ROM solutions in terms of {average quantities, e.g.,}  $L^2$-norm, and drag and lift coefficients at each timestep. 
Moreover, we {consider a pointwise (in space) quantity, i.e.,} the relative errors between the FOM and ROM velocity fields
{at} each time instance, {which is defined} as
$$
E_u(t^n) = \frac{\norm{{u}^{n} - u^n_{r_{us}}}_{L^2(\Omega)}}{\norm{{u}^{n}}_{L^2(\Omega)}},
$$ where $u^n$ is the FOM solution and $u_{r_{us}}$ represents the G-ROM, EF-ROM, VMS-EFFC-ROM, or VMS-EPFC-ROM velocity at time $t_n$. 
At the ROM level, we expand the reduced solution in the FOM space to compute these quantities.\\ 
We first compare qualitatively the velocity fields. In Figure \ref{fig:romvt1}, we plot the FOM solution (top), G-ROM (left) and EF-ROM (right) in the first row, and VMS-EFFC-ROM (left) and VMS-EPFC-ROM (right) in the second row, for $t=1$. For the differential filters we are using the same parameters as in Section \ref{sec:FOM_results}, i.e.\ $\delta_1=\delta_2=1.59\cdot 10^{-3}$, while the projection retains only the first $\overline{r}_{u}=77$ modes and filters them with $\delta=1.59\cdot 10^{-3}$. Neither G-ROM nor EF-ROM are valid options for the reconstructive task. Indeed, G-ROM presents spurious oscillations, and EF-ROM is overdiffusive. The VMS-based strategies better recover the vortex shedding behavior: the accuracy is certainly increased with respect to the other two approaches. We draw the same conclusions from Figure \ref{fig:romvt4},  which shows the FOM solution (top), G-ROM (left), and EF-ROM (right) in the first row, and VMS-EFFC-ROM (left) and VMS-EPFC-ROM (right) in the second row, for $t=4$. For this time instance, VMS-EFFC-ROM looks more accurate than the other strategies. \\
For a more quantitative analysis, in Figure \ref{fig:romnorms}, we show the squared $L^2$-norm of the reduced velocity solutions in time, to be compared with the FOM solution. On the left, we also plot the G-ROM norm, but we exclude it from the right plot to better display the performances of other strategies. The G-ROM approach yields inaccurate results. The VMS-filtered strategies look comparable, and the EF-ROM is more diffusive, as expected. However, from the right plot, we can observe that VMS-EFFC-ROM is the best approach: it is less noisy with respect to VMS-EPFC-ROM. The plot clearly shows that EF-ROM is not an accurate option in this setting.\\
We plot the drag and lift coefficients for all the approaches in Figures \ref{fig:romdrag} and \ref{fig:romdlift}, respectively. Once again, we plot the G-ROM results only in the left plots, in both figures. G-ROM fails and recovers neither drag nor lift coefficients. The performances of the other strategies can be better observed in the left plots of Figures \ref{fig:romdrag} and \ref{fig:romdlift}. From these plots, we conclude that EF-ROM is overdiffusive with constant $C_D(t)$ and $C_L(t)$. Focusing on the drag coefficient, VMS-EPFC-ROM has a very oscillatory behavior, while VMS-EFFC-ROM recovers the order of magnitude and the structure of the FOM $C_D(t)$. VMS-EFFC-ROM performs better than the other strategies for the lift coefficient: after $t=1$, we have a good agreement with the FOM lift. We can notice a slight shift with respect to the FOM, but by the end of the time interval they synchronize. However, we can observe a difference in magnitude. The VMS-EPFC-ROM shows a more inaccurate behavior: it is noisier with respect to VMS-EFFC-ROM and the shift is more pronounced. To better display these features, we present a close-up of the last two wavelengths of the lift in the right plot of Figure \ref{fig:romdlift}.

\begin{figure}[H]
  \centering
  \includegraphics[width=0.45\textwidth]{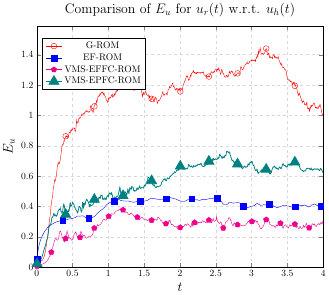}
  
  \caption{Experiment 2. Evolution of the relative error of the velocity for G-ROM, EF-ROM for $\delta=1.59 \cdot 10^{-3}$, VMS-EFFC-ROM for $\delta_1=\delta_2=1.59 \cdot 10^{-3}$, and VMS-EPFC-ROM for $\delta=1.59 \cdot 10^{-3}$ compared to the FOM solution.}
  \label{fig:romerrs}
\end{figure}

To complete the accuracy analysis, we plot the relative errors in Figure \ref{fig:romerrs}: all the DF-based strategies outperform G-ROM. VMS-EFFC-ROM is the most accurate approach, with a relative error that is six times lower than the G-ROM error, and three times lower than the VMS-EPFC-ROM error on average in time.

\begin{figure}[H]
  \centering
  \includegraphics[width=0.24\textwidth]{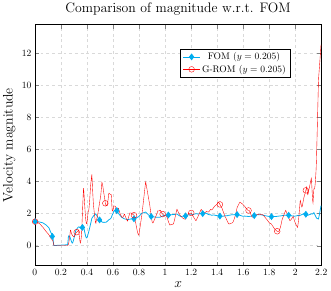}
  \includegraphics[width=0.24\textwidth]{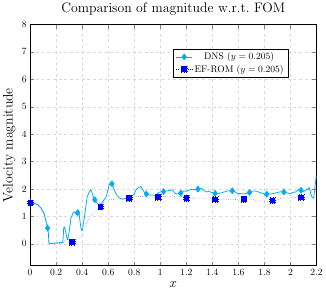}
  \includegraphics[width=0.24\textwidth]{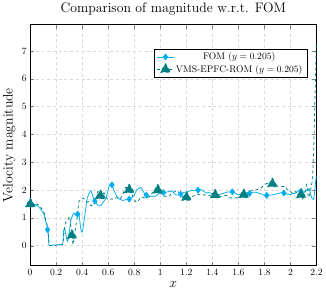}
  \includegraphics[width=0.24\textwidth]{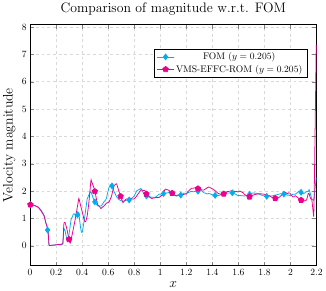}\\

  \includegraphics[width=0.24\textwidth]{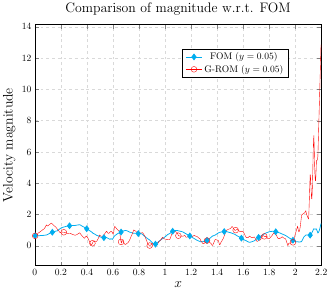}
  \includegraphics[width=0.24\textwidth]{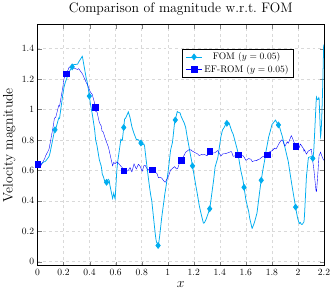}
  \includegraphics[width=0.24\textwidth]{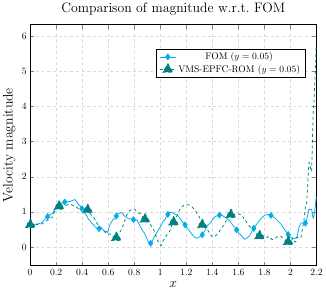}
  \includegraphics[width=0.24\textwidth]{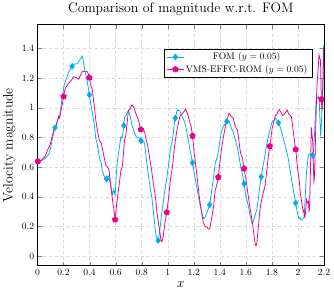}\\

   \includegraphics[width=0.24\textwidth]{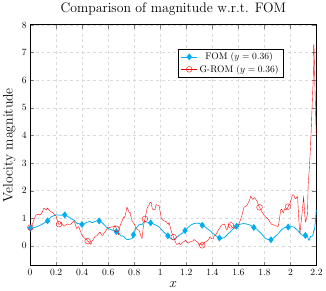}
  \includegraphics[width=0.24\textwidth]{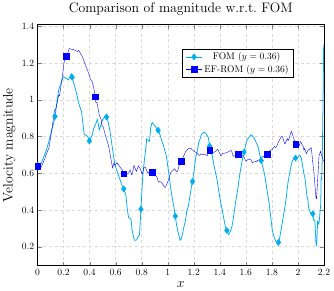}
  \includegraphics[width=0.24\textwidth]{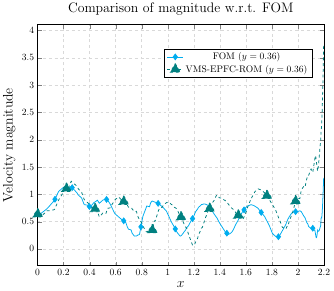}
  \includegraphics[width=0.24\textwidth]{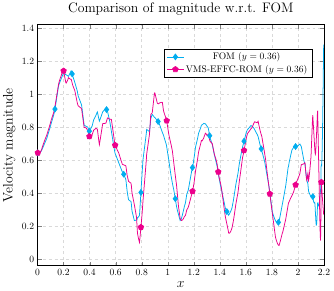}\\

  \caption{
  Experiment 2. Comparison of the ROM velocity magnitude with the FOM velocity magnitude for varying $x$ and fixed $y$. Top row from left to right: G-ROM, EF-ROM for $\delta=1.59 \cdot 10^{-3}$, VMS-EPFC-ROM for $\delta_1=\delta_2=1.59 \cdot 10^{-3}$, and VMS-EFFC-ROM for $\delta=1.59 \cdot 10^{-3}$ for $y=0.205$. Center row from left to right: G-ROM, EF-ROM for $\delta=1.59 \cdot 10^{-3}$, VMS-EPFC-ROM for $\delta_1=\delta_2=1.59 \cdot 10^{-3}$, and VMS-EFFC-ROM for $\delta=1.59 \cdot 10^{-3}$ for $y=0.05$.
  Bottom row from left to right: G-ROM, EF-ROM for $\delta=1.59 \cdot 10^{-3}$, VMS-EPFC-ROM for $\delta_1=\delta_2=1.59 \cdot 10^{-3}$, and VMS-EFFC-ROM for $\delta=1.59 \cdot 10^{-3}$ (right) for $y=0.36$.}
  \label{fig:overlineROM}
\end{figure}

To better understand the error behavior, in Figure \ref{fig:overlineROM}, we plot the ROM velocity magnitude and the FOM velocity magnitude for $t=4$ and fixed $y \in \{0.205, 0.05, 0.36\}$ for all $x \in [0,2.2]$, i.e., on three horizontal lines. First of all, G-ROM (first column of the picture) shows numerical instabilities and high oscillations at the end of the channel. These oscillations can also be observed for VMS-EPFC-ROM (third column of the picture). EF-ROM is overdiffusive and does not capture the vortex shedding. From the plot, it is clear that VMS-EFFC-ROM yields the most accurate results since it recovers the FOM solution for all the three lines. We stress, however, that VMS-EFFC-ROM also presents spurious oscillations for $y=0.205$. 
Finally, we conclude that the relative error, as in the FOM case, is not an appropriate metric for testing the accuracy of the VMS-based filters for this specific test case. Indeed, this metric is spoiled by two factors: (i) the high oscillatory behavior at the end of the channel (see, e.g., Figure \ref{fig:overlineROM}), and (ii) the shift in the lift coefficient (Figure \ref{fig:romdlift}, right plot). Nevertheless, we point out that VMS-EFFC-ROM increases the accuracy of the ROM solution for all the criteria that we considered.

\begin{remark}[On the role of $\overline{r}_u$]
\label{rem:role_proj}
For VMS-EPFC-ROM, other strategies can be chosen for projecting and defining the large scales. For example, once $\overline{r}_u$ nonzero coefficients modes are fixed, one may set to zero all the other velocity coefficients, comprising the supremizers: in this way, the supremizer are not filtered.
Another strategy discards $r'_u = r_u - \overline{r}_u$ and $r'_s = r_s - \overline{r}_s$ modes for the supremizers. In our numerical investigation, both strategies yielded worse results with respect to the approach described in Section \ref{sec:VMS-EPFC-ROM}. 
For the sake of brevity, we do not show these results.
Moreover, we tried other choices for $\overline{r}_u$, e.g., $\overline{r}_u = 
\lfloor 3r_u/4 \rfloor$ and $\overline{r}_u = 
\lfloor r_u/4 \rfloor$. In both cases, the results were worse in terms of accuracy with respect to the choice $\overline{r}_u = 
\lfloor r_u/2 \rfloor$ that we used in this section. These results are not unexpected. Indeed, for $\overline{r}_u \rightarrow r_u$, VMS-EPFC-ROM is similar to EF-ROM, while, for $\overline{r}_u \rightarrow 0$, VMS-EPFC-ROM is equivalent to the standard G-ROM, and both EF-ROM and G-ROM fail in recovering the flow structures.
\end{remark}

\begin{remark}[On the role of $\delta_1$ and $\delta_2$]
In this section, we followed the strategy of \emph{consistency} with respect to the FOM VMS-EFFC solution. Namely, we used $\delta_1=\delta_2=1.59\cdot 10^{-3}$. However, other choices can be made. For example, one can try $\delta_1 \neq \delta_2$ or 
$\delta_1=\delta_2$ but different from the FOM $\delta$ values. In both cases, our preliminary numerical investigation yielded less accurate results. Although we do not exclude that parameter-inconsistent formulations can be useful for some applications, they require a careful sensitivity analysis, which goes beyond the scope of the paper.     
\end{remark}
\begin{remark}[On the role of pressure] We note that, in our numerical investigation, the number of pressure modes $r_p$ does not affect the accuracy of the velocity fields.
We also note that $r_p = 1$ is a value that retains $99\%$ of the pressure energy. However, using only one mode for the pressure gave overdiffusive results for all the ROMs. The pressure results are shown in Appendix \ref{app:pressureROM}.
However, for brevity, we did not show the results for $r_p = 1$, and included only the ones for $r_p = 15$, which allowed us to recover the vortex shedding behavior.
\end{remark} 

Finally, we note that, since we are dealing with a nonlinear problem, the ROM depends on the FOM dimension, resulting in inefficiency during the online solution. To solve this issue, hyper-reduction techniques, such
as the empirical interpolation method (EIM), might be employed \cite{barrault2004eim}. 
However, applying these techniques was beyond the scope of this work.\\
In terms of computational time, we want to stress that the filter-based strategies are still competitive with respect to G-ROM. Indeed, the CPU time of EF and VMS-EPFC-ROM is $98\%$ of the CPU time of the G-ROM. This is due to the DF action, which helps the Newton solver converge in fewer iterations. The VMS-EFFC-ROM CPU time, however, is  $115\%$ of the G-ROM CPU time. This is because VMS-EFFC-ROM applies the DF twice, i.e., it solves two linear systems. We emphasize, however, that this slight computational overhead yields significant gains in accuracy.\

Concluding, the VMS-EFFC-ROM is the most accurate approach among the ROMs we investigated, both at the FOM and ROM levels. We also remark that both the VMS-filter strategies (i.e., VMS-EFFC-ROM and VMS-EPFC-ROM) outperform EF-ROM and G-ROM in terms of accuracy. 
\\

In Table \ref{tab:Exp}, we give an overview of the experiments of this contribution, showing their main features and the main conclusions we draw from the tests.

\begin{center}
    \begin{table}[H]
\caption{
{Summary of the numerical experiments.}}
\vspace{2mm}
\label{tab:Exp}
\centering 
\resizebox{\textwidth}{!}{
\begin{tabular}{|c|p{35mm}|p{70mm}|}
\hline
\cellcolor{gray!40}& Comparisons &Main conclusions      \\ \hline
Experiment 1 (FOM) & EF \newline VMS-EFFC \newline VMS-EPFC&  VMS-EFFC is the most accurate. \\ \hline
Experiment 2 (ROM) & G-ROM \newline EF-ROM \newline VMS-EFFC-ROM VMS-EPFC-ROM  & VMS-based strategies outperform G-ROM and EF-ROM. VMS-EFFC-ROM strategy is the most accurate. \\ \hline
\end{tabular}
}
\end{table}
\end{center}

\section{Conclusions}
\label{sec:conc}
In this work, we propose novel strategies to tackle the overdiffusivity related to the EF approach for large values of the filter radius $\delta$, which often occurs in practice. We leverage the VMS framework and propose algorithms consisting of four steps, where we (i) evolve the velocity field, (ii) decompose the evolved velocity into resolved large and small scales, (iii) filter the resolved small scales, and (iv) correct the velocity field with the filtered resolved small scales. We develop two novel strategies based on how we decompose the velocity: the VMS-EFFC, where the decomposition is performed using a differential filter, and the VMS-EPFC, which employs a postprocessing on a smaller space to define the resolved large scales. The two algorithms are numerically investigated at the FOM and ROM levels for a flow past a cylinder at Reynolds number $Re=1000$. Compared with EF and Galerkin projection approaches, the VMS-based filters yield more accurate results in terms of relative errors in time, velocity norms, pointwise errors, and flow behavior.
One of the strengths of the two new approaches is their modularity: this allows their embarrassingly easy implementation and prevents a significant increase in computational costs both offline and online.
We underline that, both at the FOM and the ROM level, VMS-EFFC is a clear winner with respect to the other strategies.\\
To the best of our knowledge, this is the first time that these approaches have been introduced both at the FOM and the ROM levels. This first investigation gives encouraging results. However, many research directions should be pursued for a deeper understanding of the topic. First, a comprehensive study on the role of the pressure and how to include it in the filtering process should be considered to increase the accuracy of the reconstruction of the pressure field. For example, one can recover the pressure with a least-squares procedure that avoids the use of supremizers, decoupling the solution for velocity and pressure, as suggested in \cite{Tomasnew}. Moreover, VMS-EPFC is sensitive to the choice of the filter radius $\delta$ and of the postprocess step, while VMS-EFFC is sensitive to the choices of $\delta_1$ and $\delta_2$. Although a comprehensive sensitivity analysis was beyond the scope of this paper but it could be useful to increase the robustness of the proposed strategies. Finally, we emphasize that this study is just a numerical investigation of the new strategies. In the future, we plan to provide numerical analysis for the VMS-based filter regularization to support the numerical findings of this manuscript.


\section*{Acknowledgments}
We acknowledge the European Union's Horizon 2020 research and innovation program under the Marie Skłodowska-Curie Actions, grant agreement 872442 (ARIA).
This study was carried out within the ``20227K44ME - Full and Reduced order modelling of coupled systems: focus on non-matching methods and automatic learning (FaReX)" project – funded by European Union – Next Generation EU  within the PRIN 2022 program (D.D. 104 - 02/02/2022 Ministero dell’Università e della Ricerca). This manuscript reflects only the authors’ views and opinions and the Ministry cannot be considered responsible for them. MS thanks the INdAM - GNCS Project ``Metodi di riduzione di modello ed approssimazioni di rango basso per problemi alto-dimensionali" (CUP\_E53C23001670001) and the ECCOMAS EYIC Grant ``CRAFT: Control and Reg-reduction in Applications for Flow Turbulence". MS further thanks the INdAM financial support as ``titolare di una borsa per l'estero dell'Istituto Nazionale di Alta Matematica". FB thanks the PRIN 2022 PNRR project ``ROMEU: Reduced Order Models for Environmental and Urban flows'' funded by the European Union -- NextGenerationEU under the National Recovery and Resilience Plan (NRRP), Mission 4 Component 2, CUP J53D23015960001. TI acknowledges support through National Science Foundation grants DMS-2012253 and CDS\&E-MSS-1953113.
The research of TC was also partially funded by the Spanish Research Agency - EU FEDER Fund Project PID2021-123153OB-C21.
The computations in this work have been performed with the RBniCS \cite{rbnics} library, which is an implementation in FEniCS \cite{fenics} of several reduced order modeling techniques; we acknowledge developers and contributors to both libraries. Computational resources were partially provided by HPC@POLITO, a project of Academic Computing within the Department of Control and Computer Engineering at the Politecnico di Torino (\href{http://hpc.polito.it}{http://hpc.polito.it}).




\bibliographystyle{abbrvurl}
\bibliography{traian.bib, maria.bib}

\appendix

\section{VMS-EFFC with Two Different Filter Radii}
\label{app:diff}
\begin{figure}[H]
 \includegraphics[width=0.49\textwidth]{Images/FOM/DNS1.png}
   \includegraphics[width=0.49\textwidth]{Images/FOM/EFFC1.png} \\
  \includegraphics[width=0.49\textwidth]{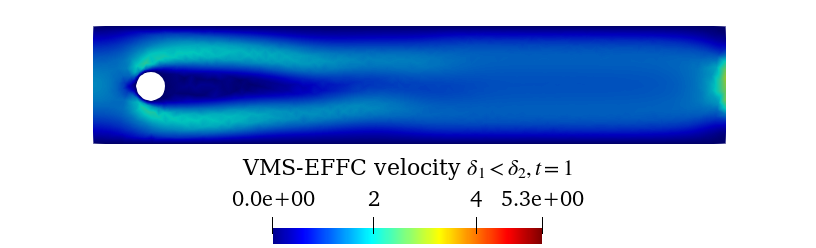}
 \includegraphics[width=0.49\textwidth]{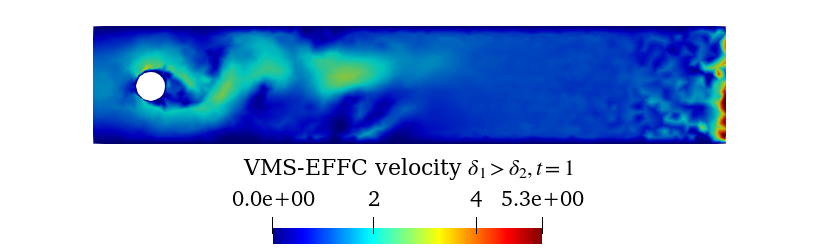}\\
 \caption{Experiment A1. Velocity profiles at $t=1$: DNS, VMS-EFFC for $\delta_1=\delta_2=1.59 \cdot 10^{-3}$, VMS-EFFC for $\delta_1=1.59 \cdot 10^{-3}$ and $\delta_2 = 10 \cdot \delta_1$, and VMS-EFFC for $\delta_1=1.59 \cdot 10^{-3}$ and $\delta_2 = 10^{-1} \delta_1$, top left and right, bottom left, and right plots, respectively.}
 \label{fig:diffvt1a1}
\end{figure}

\begin{figure}[H]
 \includegraphics[width=0.49\textwidth]{Images/FOM/DNS4.png}
   \includegraphics[width=0.49\textwidth]{Images/FOM/EFFC4.png} \\
  \includegraphics[width=0.49\textwidth]{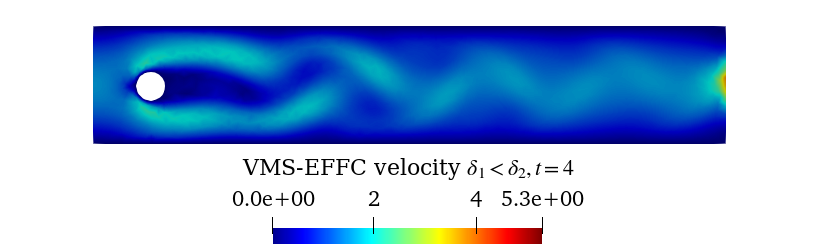}
 \includegraphics[width=0.49\textwidth]{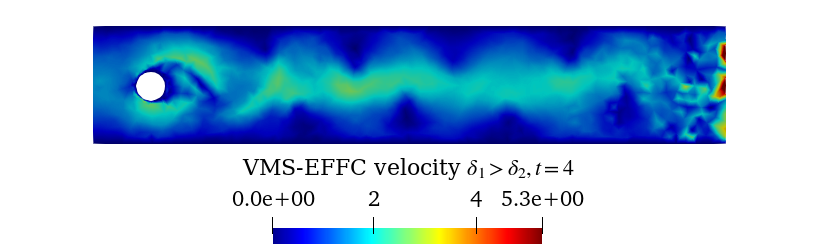}\\
 \caption{Experiment A1. Velocity profiles at $t=4$: DNS, VMS-EFFC for $\delta_1=\delta_2=1.59 \cdot 10^{-3}$, VMS-EFFC for $\delta_1=1.59 \cdot 10^{-3}$ and $\delta_2 = 10 \cdot \delta_1$, and VMS-EFFC for $\delta_1=1.59 \cdot 10^{-3}$ and $\delta_2 = 10^{-1} \delta_1$, top left and right, bottom left, and right plots, respectively.}
 \label{fig:diffvt4a1}
\end{figure}

We perform a numerical investigation of the impact of using two different filter radii, $\delta_1 \neq \delta_2$, for the VMS-EFFC approach. \\
The numerical setting is described in Section \ref{sec:gen_setting}. We propose two different test cases: (i) with $\delta_1 = 1.59 \cdot 10^{-3}$ and $\delta_2 = 1.59 \cdot 10^{-2}$, i.e., $\delta_1 < \delta_2$, with a one order of magnitude difference, and (ii) with $\delta_1 = 1.59 \cdot 10^{-3}$ and $\delta_2 = 1.59 \cdot 10^{-4}$, i.e., $\delta_1 > \delta_2$, with a one order of magnitude difference.\\
In Figures \ref{fig:diffvt1a1} and \ref{fig:diffvt4a1}, we show the DNS velocity (top left plot) and we compare it to the VMS-EFFC strategy for $\delta_1 = \delta_2$ (top right plot), for $\delta_1 < \delta_2$ (bottom left plot) and for $\delta_1 > \delta_2$ (bottom right plot), for $t=1$ and $t=4$, respectively. It is clear that filtering the small scale with a larger $\delta_2$ can produce overdiffusive results, while using a smaller $\delta_2$ can be insufficient to alleviate the numerical oscillations of the flow. For the sake of completeness, we also include the case for $\delta_1 = \delta_2$, which was already analyzed in Section \ref{sec:FOM_results}.

\begin{figure}[H]
 \centering
 \includegraphics[width=0.472\textwidth]{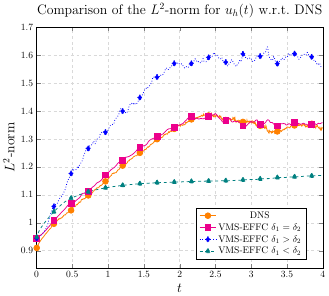}
 \caption{Experiment A1. Evolution of the $L^2$-norms of the velocity for DNS for $\delta=1.59 \cdot 10^{-3}$, VMS-EFFC for $\delta=1.59 \cdot 10^{-3}$ and $\delta_2 = 1.59 \cdot 10^{-2}$, and VMS-EFFC for $\delta=1.59 \cdot 10^{-3}$ and $\delta_2 = 1.59 \cdot 10^{-4}$.}
 \label{fig:normsa1}
\end{figure}

\begin{figure}[H]
 \centering
 \includegraphics[width=0.49\textwidth]{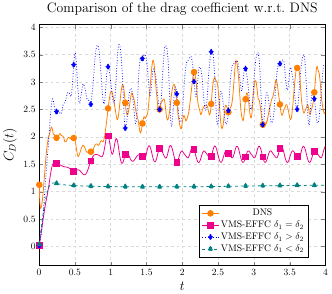}
 \includegraphics[width=0.49\textwidth]{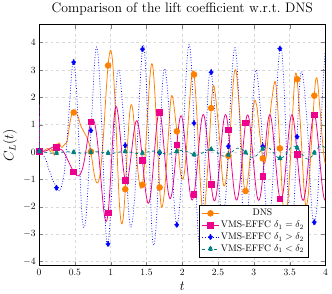}
 \caption{Experiment 1. $C_D(t)$ and $C_L(t)$ in time for DNS for $\delta=1.59 \cdot 10^{-3}$, VMS-EFFC for $\delta=1.59 \cdot 10^{-3}$ and $\delta_2 = 1.59 \cdot 10^{-2}$, and VMS-EFFC for $\delta=1.59 \cdot 10^{-3}$ and $\delta_2 = 1.59 \cdot 10^{-4}$, left and right plots, respectively.}
 \label{fig:draglifta1}
\end{figure}

In this numerical setting, using different values for $\delta_1$ and $\delta_2$ does not increase the accuracy. This is the conclusion we draw from Figure \ref{fig:normsa1}. The left plot shows that the best strategy to recover the $L^2$-norm of the DNS simulation is to use $\delta_1 = \delta_2$. The overdiffusive action is also illustrated by Figure \ref{fig:draglifta1}: The use of $\delta_1 < \delta_2$ damps out the oscillations of the drag and lift coefficients. 

Using $\delta_1 > \delta_2$, helps to recover the magnitude of the drag and lift coefficients. However, $\delta_1 > \delta_2$, yields significantly larger values of energy norms compared to the DNS, and it gives a worse representation of the velocity field in terms of the relative error, too. 
From this analysis, it is clear that the choice of using $\delta_1 = \delta_2$ leads to a better qualitative performance, including a better norm reconstruction. We stress that this conclusion holds only for our numerical setting, and a proper choice of $\delta_1$ and $\delta_2$ may increase the accuracy of the velocity field reconstruction. However, finding the optimal values for $\delta_1$ and $\delta_2$ was beyond the scope of the manuscript. The pressure results for this experiment are included in Appendix \ref{app:pressure}.

\section{Additional FOM Results: Pressure}
\label{app:pressure}
\begin{figure}[H]
  \includegraphics[width=0.49\textwidth]{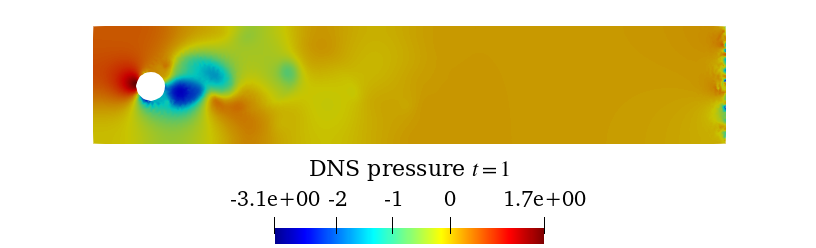}
      \includegraphics[width=0.49\textwidth]{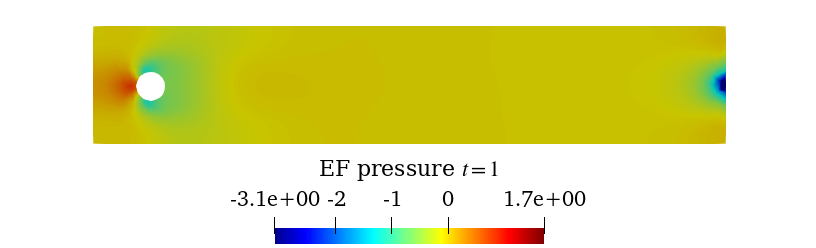} \\
    \includegraphics[width=0.49\textwidth]{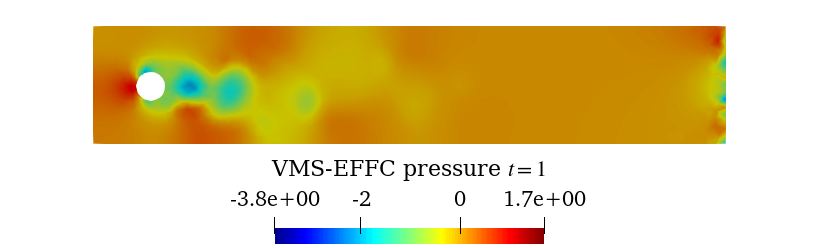}
  \includegraphics[width=0.49\textwidth]{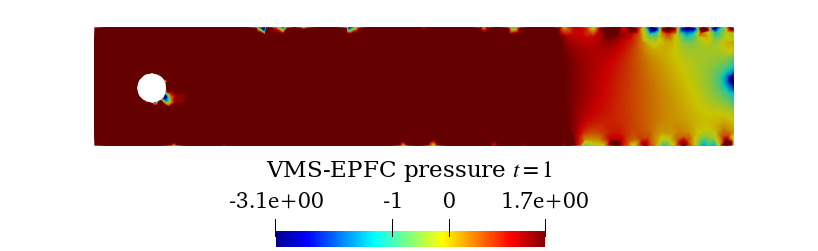}\\
 \caption{Experiment 1. Pressure profiles at $t=1$: DNS, EF for $\delta=1.59 \cdot 10^{-3}$, VMS-EFFC for $\delta_1=\delta_2=1.59 \cdot 10^{-3}$, and VMS-EPFC for $\delta=1.59 \cdot 10^{-3}$ and $\overline{r}_u = 15$, from left to right, and top to bottom.}
  \label{fig:pt1}
\end{figure}

\begin{figure}[H]
  \includegraphics[width=0.49\textwidth]{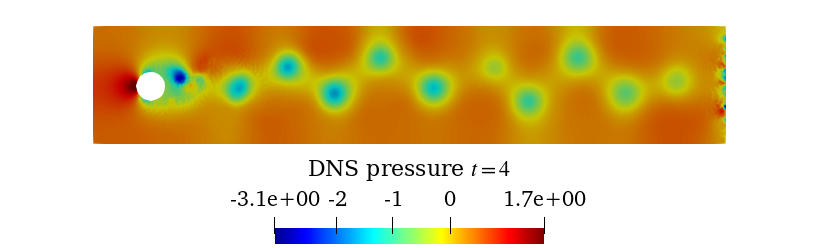}
      \includegraphics[width=0.49\textwidth]{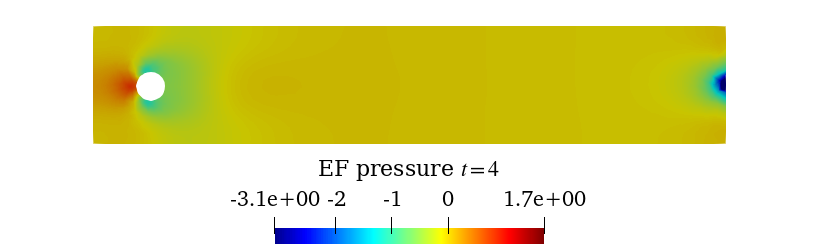} \\
    \includegraphics[width=0.49\textwidth]{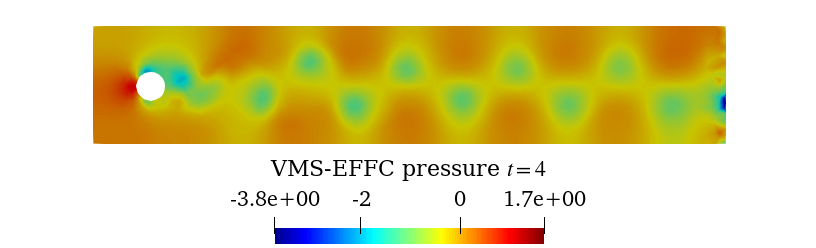}
  \includegraphics[width=0.49\textwidth]{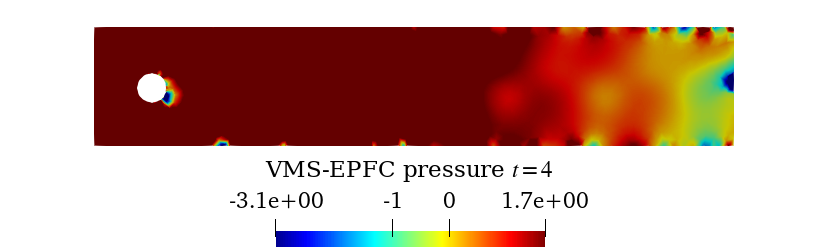}\\
 \caption{Experiment 1. Pressure profiles at $t=4$: DNS, EF for $\delta=1.59 \cdot 10^{-3}$, VMS-EFFC for $\delta_1=\delta_2=1.59 \cdot 10^{-3}$, and VMS-EPFC for $\delta=1.59 \cdot 10^{-3}$, from left to right, and top to bottom.}
  \label{fig:pt4}
\end{figure}

\begin{figure}[H]
  \centering
  \includegraphics[width=0.39\textwidth]{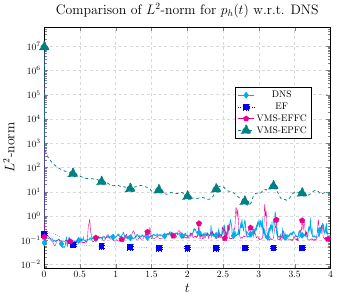}
  \caption{Experiment 1. Evolution of the $L^2$-norms of the pressure for DNS, EF for $\delta=1.59 \cdot 10^{-3}$, VMS-EFFC for $\delta_1=\delta_2=1.59 \cdot 10^{-3}$, and VMS-EPFC for $\delta=1.59 \cdot 10^{-3}$.}
  \label{fig:normsp}
\end{figure}

\begin{figure}[H]
  \includegraphics[width=0.49\textwidth]{Images/FOM/DNSp1.png}
      \includegraphics[width=0.49\textwidth]{Images/FOM/EFFCp1.png} \\
    \includegraphics[width=0.49\textwidth]{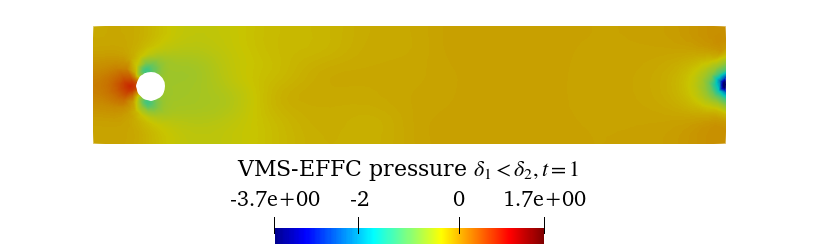}
  \includegraphics[width=0.49\textwidth]{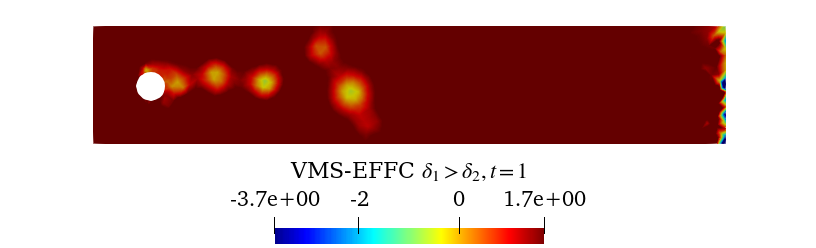}\\
 \caption{Experiment A1. Pressures at $t=1$: DNS, VMS-EFFC for $\delta_1=\delta_2=1.59 \cdot 10^{-3}$, VMS-EFFC for $\delta_1=1.59 \cdot 10^{-3}$ and $\delta_2 = 10 \cdot \delta_1$, VMS-EFFC for $\delta_1=1.59 \cdot 10^{-3}$ and $\delta_2 = 10^{-1}\cdot \delta_1$, from left to right, and top to bottom.}
  \label{fig:pt1a1}
\end{figure}

\begin{figure}[H]
  \includegraphics[width=0.49\textwidth]{Images/FOM/DNSp4.png}
      \includegraphics[width=0.49\textwidth]{Images/FOM/EFFCp4.png} \\
    \includegraphics[width=0.49\textwidth]{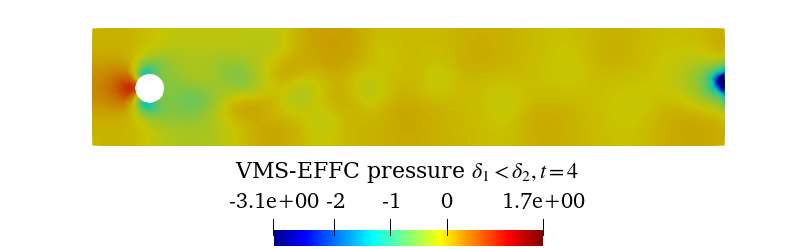}
  \includegraphics[width=0.49\textwidth]{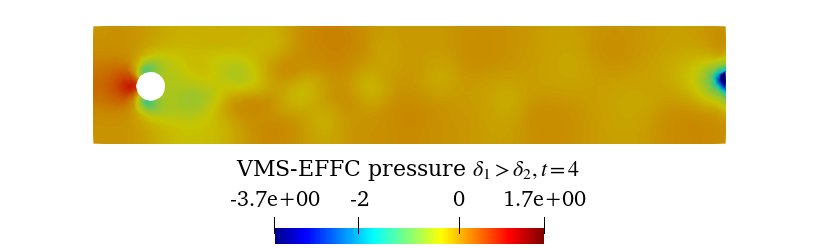}\\
 \caption{Experiment A1. Pressures at $t=4$: DNS, VMS-EFFC for $\delta_1=\delta_2=1.59 \cdot 10^{-3}$, VMS-EFFC for $\delta_1=1.59 \cdot 10^{-3}$ and $\delta_2 = 10 \cdot \delta_1$, VMS-EFFC for $\delta_1=1.59 \cdot 10^{-3}$ and $\delta_2 = 10^{-1}\cdot \delta_1$, from left to right, and top to bottom.}
  \label{fig:pt4a1}
\end{figure}

\begin{figure}[H]
  \centering
  \includegraphics[width=0.375\textwidth]{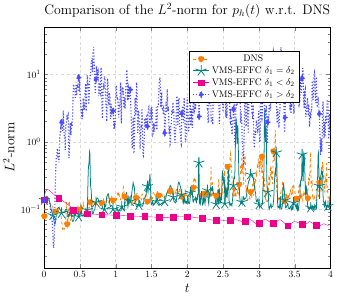}
  \caption{Experiment A1. Evolution of the squared $L^2$-norms of the pressure DNS, EF for $\delta=1.59 \cdot 10^{-3}$, VMS-EFFC for $\delta=1.59 \cdot 10^{-3}$ and $\delta_2 = 1.59 \cdot 10^{-2}$, and VMS-EFFC for $\delta=1.59 \cdot 10^{-3}$ and $\delta_2 = 1.59 \cdot 10^{-4}$, left and right plots, respectively.}
  \label{fig:norms_pa1}
\end{figure}

This appendix collects the FOM results of the pressure fields corresponding to Experiment 1 (Figures \ref{fig:pt1}, \ref{fig:pt4}, \ref{fig:normsp}) and Experiment A1 (Figures \ref{fig:norms_pa1}, \ref{fig:pt1a1}, \ref{fig:pt4a1}).


The results are consistent with the behavior of the velocity fields in the respective experiments. In what follows, we summarize the main features of the pressure results.

\begin{itemize}
    \item[$\circ$] Experiment 1. The best approach is VMS-EFFC, which is capable of representing the DNS pressure field. In contrast, EF and VMS-EPFC fail to represent the DNS pressure field. This is clearly supported by the qualitative plots in Figures \ref{fig:pt1} and \ref{fig:pt4}, and the energy norm in the left plot of Figure \ref{fig:normsp}.
    \item[$\circ$] Experiment A1. From the qualitative plots of Figures \ref{fig:pt1a1} and \ref{fig:pt4a1}, and the energy norm plot on the left of Figure \ref{fig:norms_pa1}, it is clear that VMS-EFFC with $\delta_1 = \delta_2$ is the best approach for pressure representation. Indeed, it recovers the DNS pressure behavior and magnitude of the $L^2$-norm, while the strategies that use $\delta_1 \neq \delta_2$. 
\end{itemize}


\section{Additional ROM Results: Pressure}
\label{app:pressureROM}
In this appendix, we present the results of the ROM reconstruction of the pressure for Experiment 2. These results are consistent with the velocity results of Section \ref{sec:ROM_results}. 
We briefly comment on these results, emphasizing the salient points. First, we note that all the strategies fail to represent the pressure fields. 
\begin{figure}[H]
\centering
  \includegraphics[width=0.49\textwidth]{Images/FOM/EFFCp1.png}\\
      \includegraphics[width=0.49\textwidth]{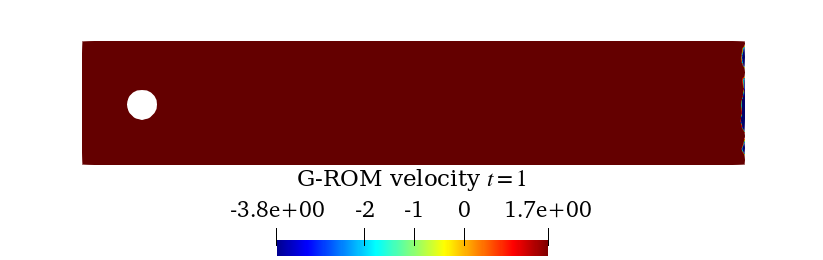}
    \includegraphics[width=0.49\textwidth]{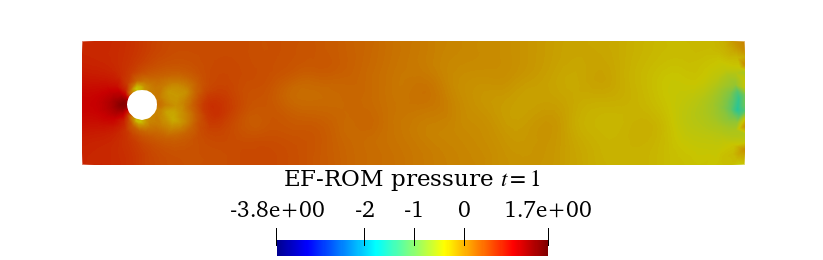}
  \includegraphics[width=0.49\textwidth]{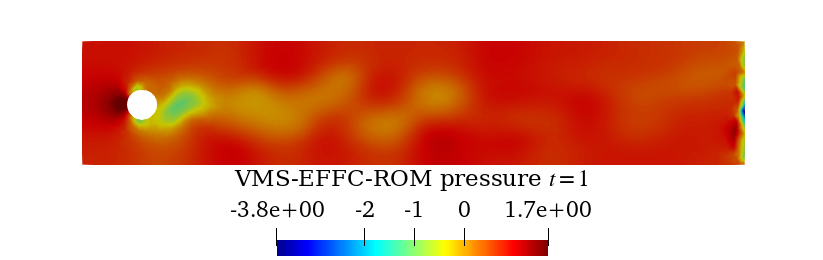}   \includegraphics[width=0.49\textwidth]{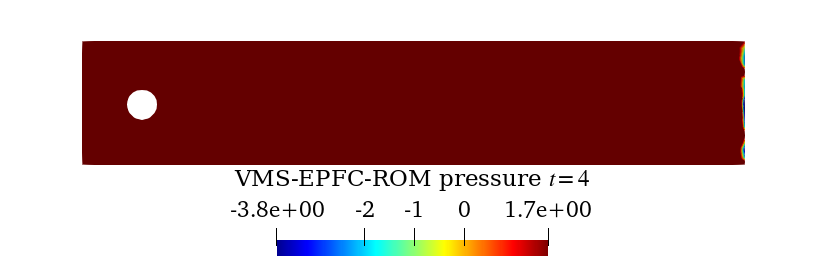}\\
  \caption{Experiment 2. Pressure profiles at $t=1$: FOM (VMS-EFFC) for $\delta=1.59 \cdot 10^{-3}$, G-ROM, EF-ROM for $\delta=1.59 \cdot 10^{-3}$, VMS-EFFC-ROM $\delta_1=\delta_2=1.59 \cdot 10^{-3}$, and VMS-EPFC-ROM for $\delta=1.59 \cdot 10^{-3}$ and $\overline{r}_u = 77$, top, middle left, middle right, bottom left, and bottom right plots, respectively.}
  \label{fig:romvt1p}
\end{figure}
\begin{figure}[H] \centering
  \includegraphics[width=0.49\textwidth]{Images/FOM/EFFCp4.png}\\
      \includegraphics[width=0.49\textwidth]{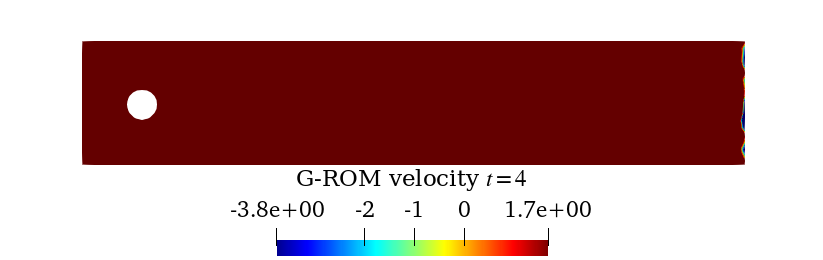}
    \includegraphics[width=0.49\textwidth]{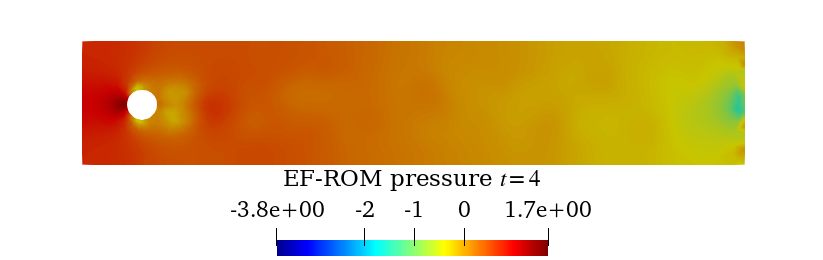}
  \includegraphics[width=0.49\textwidth]{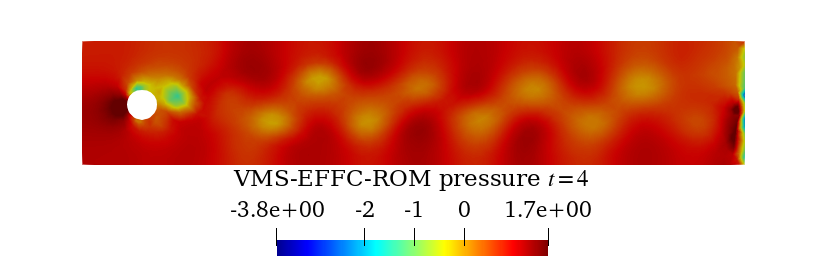}   \includegraphics[width=0.49\textwidth]{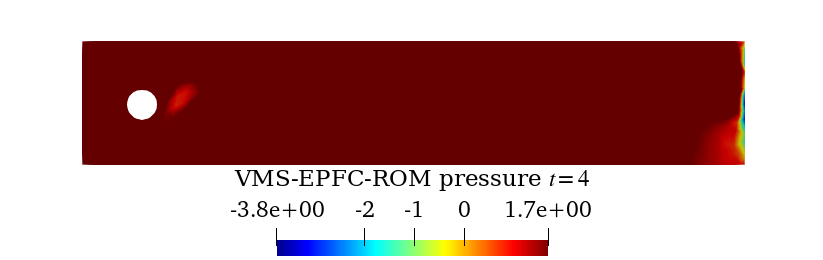}\\
  \caption{Experiment 2. Pressure profiles at $t=4$: FOM (VMS-EFFC) for $\delta=1.59 \cdot 10^{-3}$, G-ROM, EF-ROM for $\delta=1.59 \cdot 10^{-3}$, VMS-EFFC-ROM $\delta_1=\delta_2=1.59 \cdot 10^{-3}$, and VMS-EPFC-ROM for $\delta=1.59 \cdot 10^{-3}$ and $\overline{r}_u = 77$, top, middle left, middle right, bottom left, and bottom right plots, respectively.}
  \label{fig:romvt4p}
\end{figure}

\begin{figure}[H]
  \centering
  \includegraphics[width=0.44\textwidth]{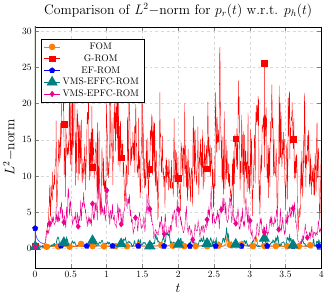}
  \includegraphics[width=0.44\textwidth]{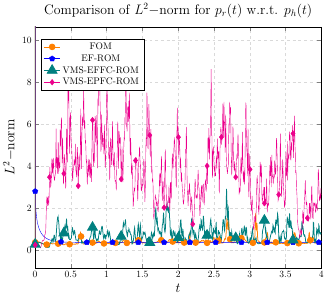}
  
  \caption{Experiment 2. \emph{Left and right}: Evolution of the squared $L^2$-norms of the pressure for G-ROM, EF-ROM for $\delta=1.59 \cdot 10^{-3}$, VMS-EFFC for $\delta_1=\delta_2=1.59 \cdot 10^{-3}$, VMS-EPFC for $\delta=1.59 \cdot 10^{-3}$ and $\overline{r}_{u}=77$, and the FOM solution (VMS-EFFC) for $\delta=1.59 \cdot 10^{-3}$. The left plot shows G-ROM, which is omitted in the central plot for clarity}
  \label{fig:romnormsP}
\end{figure}
This is the conclusion we draw from Figures \ref{fig:romvt1p} and \ref{fig:romvt4p}, where we compare the ROM solutions of the four different strategies with the FOM pressure solution for $t=1$ and $t=4$, respectively. We note that these plots are scaled with respect to the FOM solution magnitude. These plots show that all the ROM techniques fail to recover the magnitude of the pressure. 
We note, however, that the VMS-EFFC-ROM strategy recovers the vortex shedding behavior.
Quantitatively speaking, there is no accurate strategy for this computational setting. We can observe, however, that VMS-EFFC-ROM reaches the minimum relative error values, as expected.
These results suggest that the VMS-filter strategies are improvable in terms of recovering the pressure. Although this investigation was beyond the scope of this paper, we plan to pursue it in future work, as mentioned in Section \ref{sec:conc}.

\end{document}